\documentclass[12pt, a4paper]{article}

\usepackage{amsmath,amsthm,amsfonts,amssymb,latexsym,stmaryrd}

\begin{document}

\newtheorem*{theo}{Theorem}
\newtheorem*{pro} {Proposition}
\newtheorem*{cor} {Corollary}
\newtheorem*{lem} {Lemma}
\newtheorem{theorem}{Theorem}[section]
\newtheorem{corollary}[theorem]{Corollary}
\newtheorem{lemma}[theorem]{Lemma}
\newtheorem{proposition}[theorem]{Proposition}
\newtheorem{conjecture}[theorem]{Conjecture}
 \newtheorem{definition}[theorem]{Definition}
 \newtheorem{remark}[theorem]{Remark}
 \newtheorem{example}[theorem]{Example}
\newcommand{\Naturali}{{\mathbb{N}}}
\newcommand{\Reali}{{\mathbb{R}}}
\newcommand{\Complessi}{{\mathbb{C}}}
\newcommand{\Toro}{{\mathbb{T}}}
\newcommand{\Relativi}{{\mathbb{Z}}}
\newcommand{\HH}{\mathfrak H}
\newcommand{\KK}{\mathfrak K}
\newcommand{\LL}{\mathfrak L}
\newcommand{\as}{\ast_{\sigma}}
\newcommand{\tn}{\vert\hspace{-.3mm}\vert\hspace{-.3mm}\vert}
\def\A{{\cal A}}
\def\B{{\cal B}}
\def\H{{\cal H}}
\def\K{{\cal K}}
\def\L{{\cal L}}
\def\N{{\cal N}}
\def\M{{\cal M}}
\def\O{{\cal O}}
\def\P{{\cal P}}
\def\S{{\cal S}}
\def\T{{\cal T}}
\def\U{{\cal U}}

\title{On twisted Fourier analysis and convergence \\
of Fourier series on discrete groups
}

\author{Erik B\'edos$^*$,
Roberto Conti$^{**}$ \\}
\date{February 21, 2008}
\maketitle
\markboth{}{}

\renewcommand{\sectionmark}[1]{}
\begin{abstract} 
We study norm convergence and summability
of Fourier series in the setting of reduced twisted group $C^*$-algebras of discrete groups. 
For amenable groups,  F\o lner nets  give the key to Fej\'er summation. 
We show that Abel-Poisson summation holds for a large class of groups,
including e.g.\ all Coxeter groups and all Gromov hyperbolic groups. As a tool in our presentation,
 we introduce notions of polynomial and subexponential H-growth for
 countable groups w.r.t.\ proper scale functions, usually chosen as length functions.
 These coincide with the classical notions of growth in the case of amenable groups.


\end{abstract}

\vfill
\thanks{\noindent $^*$  partially supported by the Norwegian Research
Council.\par
 \noindent Address: Institute of Mathematics, University of
Oslo,
P.B. 1053 Blindern,\\ 0316 Oslo, Norway. E-mail: bedos@math.uio.no. \\

\noindent $^{**}$
Address:  Mathematics, School of Mathematical and Physical  Sciences, 
University of Newcastle, Callaghan,  NSW, 2308, Australia. \\
E-mail: roberto.conti@newcastle.edu.au.\par}

\newpage

\bigskip

\section{Introduction} 

Let $T$ be a compact abelian group and $G=\widehat{T}$ denote its dual 
Letting $C(T)$, the continuous complex functions on $T$,
act as multiplication operators on $L^2(T)$, and identifying $L^2(T)$ with $\ell^2(G)$ via Fourier transform, 
one obtains $C^*_{r}(G)$, the reduced group C$^*$- algebra of $G$, which is generated by the translation operators  
$\lambda (g)$ on $\ell^2(G)$. 
In the same way, $L^\infty(T)$ corresponds to $vN(G)$, the group von Neumann algebra of $G$. In this picture, 
the uniform norm $\| \cdot \|_{\infty}$ becomes the operator norm $\| \cdot \|$.
Now, $C^*_{r}(G)$ and $vN(G)$ make sense for any discrete group $G$ and 
may then be thought of as  dual objects associated with $G$. Ever since the pioneering work of Murray and von Neumann, 
such group algebras (and their locally compact analogues) have been an important source of examples in operator algebra theory. 
More recently, they have also inspired several concepts, results and conjectures in noncommutative geometry, 
as illustrated in \cite{Co, V, GuHi}. 
In many situations (see e.g.\ \cite{Bel1, Bel2, CHMC,  MM1, GL, MM2, Mat, Lu1, Lu2}), it
appears to be useful to
consider also twisted versions of these algebras, $C_{r}^*(G, \sigma)$ and $vN (G, \sigma)$, where  $\sigma$ is a 2-cocycle
on $G$ with values in the unit circle $\Toro$, the generators being now twisted translation operators $\Lambda_{\sigma}(g)$
acting on $\ell^2(G)$ and satisfying $\Lambda_{\sigma}(g)\Lambda_{\sigma}(h)=\sigma(g,h)\Lambda_{\sigma}(gh)$ for all
$g, h \in G$. Except in trivial cases, these
twisted algebras are noncommutative, even if $G=\widehat{T}$ is abelian, and can not be studied by classical methods. 
A simple, but very popular example is $G=\Relativi^N$ with $\sigma_{\Theta}: \Relativi^N \times \Relativi^N \to \Toro$ 
given by $\sigma_{\Theta}(x, y)=\text{exp}(i\, x^t \, \Theta \, y) $ for some $N\times N$ real matrix $\Theta$, the resulting
$C_{r}^*(G, \sigma_{\Theta})$ being then called a noncommutative $N$-torus whenever $\sigma_{\Theta}$ is not symmetric. 

\medskip 
Given any result in harmonic analysis on a compact abelian group $T$ which may be reformulated as a result about $C_{r}^*(\widehat{T})$ 
(or $vN(\widehat{T})$), one may wonder whether this result carries over to $C_{r}^*(G, \sigma)$
(or $vN(G, \sigma)$). As our starting point, we consider the  basic problem 
that the Fourier series of a function $f \in C(T)$ does not necessarily converge uniformly (to $f$).
The usual way to remedy for this defect,
at least when $T=\Toro,$  
 is either to  assume that $f$ is $C^1$, or more generally, that $\widehat{f} \in \ell^1(\Relativi)$,
 or to follow ideas of Abel, Cesaro, Poisson and Fej\'er, and introduce 
 different kind of summation processes for Fourier series. 
  
  \medskip Let us briefly review here these summation processes when $T=\Toro$, so $G=\Relativi$. 
For each $k \in \Relativi$ set $e_{k}(z)= z^k \  (z \in \Toro)$, so that for $f \in C(\Toro)$, 
the (formal) Fourier series of $f$ is given by $\sum_{k \in \Relativi}\widehat{f}(k)e_{k}.$
  
\medskip Let $(\varphi_{n})$ be a sequence in $\ell^2(\Relativi)$. For each $n \in \Naturali$, set
$$M_{n}(f) := \sum_{k \in \Relativi}\varphi_{n}(k)\widehat{f}(k)e_{k}, \quad f \in C(\Toro),$$
this series being clearly absolutely convergent with respect  to $\| \cdot \|_{\infty}$.
Each $M_{n} $ is then a bounded linear map from $(C(\Toro),  \| \cdot \|_{\infty})$ into itself, satisfying  $\| M_{n} \| \leq  \| \varphi_{n} \|_{2}.$  
     
     Elementary  analysis shows that  $M_{n}(f)$ converges uniformly (necessarily to $f$) for all $f \in C(\Toro)$ 
      if and only if  $\varphi_{n} \to 1$ pointwise on $\Relativi$
      and   $\sup_{n}\| M_{n} \|$ is finite,
      in which case one could say that $C(\Toro)$ has the summation property with respect to $(\varphi_{n}).$
      The main difficulty in a concrete situation is to compute  
      the operator norms $\| M_{n} \|,$ or at least to get good estimates for them (the problem being that 
      in cases of interest $ \| \varphi_{n} \|_{2} \to \infty$). 
      
      The usual convergence problem of Fourier series consists of looking at
      $\varphi_{n}(k)=d_{n}(k):= 1$ if $| k | \leq n$ and $0$ otherwise, 
      in which case
 $\| M_{n} \|  \to \infty $. 
      In the case of Fej\'er summation, one considers instead
      $\varphi_{n}(k)=f_{n}(k):= 1 - \frac {| k | } {n}$ if $| k | \leq n-1$ and $0$ otherwise.
      Then  $\| M_{n} \| = 1$ for all $n$, and it follows that the Fourier series of any $f$ in $C(\Toro)$ is
       uniformly Fej\'er summable to $f.$
      For Abel-Poisson summation, one picks a sequence $(r_{n})$ in the interval $(0,1)$ converging to $1$
       and considers  $\varphi_{n}(k)=p_{n}(k):= r_{n}^{| k |}.$ 
     Alternatively, one can introduce $p_{r}(k)= r^{| k |}$ for  $r \in (0,1)$ and the associated operator $M_{r}$
     defined in the obvious way, and let $r \to 1$ as is usually done. 
     (We will use nets instead of sequences in the sequel to accomodate for such situations). 
      Again $\| M_{n} \| = 1 \, = \| M_{r} \| $ 
      for all $n$ and all $r$, hence the Fourier series of any $f$ in $C(\Toro)$ is uniformly Abel-Poisson summable to $f.$   
  
   The proofs of these results usually invoke the  Fej\'er kernel $F_n$ (writing $M_{n}(f)=F_{n}*f$) and  
   the Poisson kernel $P_r$ (writing $M_{r}(f)=P_{r}*f$), 
   which have nicer behaviour than the Dirichlet kernel $D_n$. In fact, as they are both non-negative,
their Fourier transforms 
   $\widehat{F_{n}} = f_{n}$ and $\widehat{P_{r}}= p_{r}$  are {\it positive definite} functions on $\Relativi$ 
   (while  $\widehat{D_n} = d_{n}$ is not), and this sole fact implies that $\| M_{n}\|=1$ for all $n$.  
   This way of establishing the key property for a summation process is not new (and is also relevant when considering
   $L^1(\Toro)$ instead of $C(\Toro)$): see e.g.\ \cite[Section 1]{ER}. It has the interesting feature that it generalizes to
   a broader context.

 \medskip Consider now some reduced twisted group $ C^*$-algebra $A= C^*_{r}(G, \sigma)$ associated to a discrete group $G$.
 To each element $x \in A$ on may naturally define its (formal) Fourier series $\sum_{g\in G}\widehat{x}(g)\Lambda_{\sigma}(g)$ 
 (see Section 2). In this paper, we address the following problems:
 
   \smallskip i) giving conditions ensuring that this series converges in operator norm (necessarily to $x$).
 
 \smallskip ii)  establishing the existence of summation processes on $A$.

 \smallskip  Concerning i), it is clear that the condition $\widehat{x}\in \ell^1(G)$  provides one natural answer. In classical
  theory, the degree of smoothness of $x$ is reflected in some stronger decay condition on $\widehat{x}$, with $x \in C^\infty$ corresponding
  to $\widehat{x}$ having rapid decay. Inspired by the work of
    P. Jolissaint \cite{Jol} on groups with the rapid decay property (with respect to some length function) and its twisted
    version \cite{Cha}, we illustrate in Section 3 how 
   i) may also be answered by introducing some decay conditions (w.r.t.\ some weight function on $G$), which involve not only the 
   Fourier coefficients of an element, but also $G$. 
To illustrate the applicability of these conditions, we use ideas from U.\ Haagerup's  
   paper \cite{Haa1} and introduce notions of
   polynomial and subexponential H-growth  for a countable group $G$. These notions of H-growth are defined w.r.t.\ to a proper 
   (scale) function on $G$, which is commonly taken to be a length function. Instead of using (the square root of) the cardinality $|E|$
   of a finite nonempty subset $E$ of $G$ to measure growth, we use  the "Haagerup content" $c(E)$ of $E$, which may be defined as follows :
   $$c(E)= \sup \Big\{ \,Ê\Big\|Ê\sum_{g\in E} a_{g} \lambda(g) \Big\| \; |Ê\;  \; a_{g} \in \Complessi \; (g \in E), \;
    \sum_{g \in E} |a_{g}|^2 = 1 \Big\}$$
    i.e.\ $c(E)$ is the norm of the natural embedding of $\ell^2(E)$ into $C_{r}^*(G).$
    This number satisfies $1 \leq c(E)\leq |E|^{1/2}$.
   Moreover,  $c(E)=|E|^{1/2}$ for all $E$ whenever $G$ is amenable \cite{Pa, Pi}. Hence, one recovers
 the usual notions of growth  in the case of amenable groups (and  proper length functions).   
   Our main contribution concerning problem i) may then be summarized as follows: 
  
  \begin{theorem} \label{I1} Let $L: G\to [0, \infty)$ be a proper function. 
   
   \smallskip If $G$ has polynomial H-growth (w.r.t.\ $L$), then there exists some $s >0 $ such that the Fourier series of 
   $x \in C_{r}^*(G,\sigma)$ converges  to $x$
   in operator norm  whenever $\sum_{g\in G} | \widehat{x}(g)|^2 \, (1 + L(g))^{s} < \infty .$
   
\smallskip If $G$ has subexponential H-growth (w.r.t.\ $L$), then the Fourier series of $x \in C_{r}^*(G,\sigma)$ converges  to $x$
   in operator norm  whenever there exists some $t > 0$ such that  $\sum_{g\in G} |\widehat{x}(g)|^2 \,\exp(t L(g)) < \infty .$
   \end{theorem} 
  To mention just one example here, the first statement in this theorem applies when $G$ is a free group on $n$ generators and $L$ is the canonical word-length
  function on $G$, in which case one may choose any $s > 3$. 
   
  \bigskip 
  As an intermediate step before discussing problem ii), we study multipliers on $C_{r}^*(G,\sigma)$ in Section 4, and pay attention
  to those which transform each element of $C_{r}^*(G,\sigma)$ into an element having an operator norm convergent Fourier series. These
special multipliers are used to define summation processes in Section 5. It has  been known already   
   since the work of G.\ Zeller-Meier  that some analogue of Fej\'er summation for Fourier series exists
    when the group is amenable (cf.\ \cite[Proposition 5.8]{ZM}; see also \cite{Exe} for the untwisted case and \cite{Wea} for $G=\Relativi^2$ with a twist). 
The direct analogue of Fej\'er summation may in fact be obtained in this case after picking a F\o lner net  in $G$,
  the existence of such a net being equivalent to the amenability of the group. The precise statement is as follows.
   
  \begin{theorem} \label{I2} Assume that $\{F_\alpha\}$ is  a F\o lner net of finite subsets for $G.$  
Then 
$$  \sum_{g \in  F_\alpha \cdot F_\alpha^{-1}} 
\frac{ |g F_\alpha  \cap F_\alpha|}{|F_\alpha|}  \, \,Ê\widehat{x}(g) \Lambda_\sigma(g) 
\underset\alpha \longrightarrow x \ \ \text{(in operator norm)}$$
for all $x \in C^*_{r}(G,\sigma)$. 
\end{theorem}
If $G = \Relativi$ and $\sigma =1$,  then choosing $F_n = \{0,1,\ldots,n-1\}$ 
just gives the classical Fej\'er summation theorem.

   \medskip 
   The  analogue of Abel-Poisson summation is more troublesome, unless 
   the group is $\Relativi^N$ for some $N \in \Naturali$. In this case, we
   show the following : 
 
  \begin{theorem}  \label{I3} Let $p$ be 1 or 2, and $| \cdot |$ denote the usual $p$-norm on $\Relativi^N$. 
   Let $r \in (0,1)$. Then
$$  \sum_{g \in \Relativi^N} 
  r^{| g |}\,Ê\widehat{x}(g) \Lambda_\sigma(g) 
\underset{r\to 1^{-}} \longrightarrow x \; \; \text{(in operator norm)}$$ 
for all $x \in C^*_r(\Relativi^N ,\sigma)$.
\end{theorem}  
 
  Actually, we also show that the twisted analogue of the Abel-Poisson summation theorem 
  holds for a large class of groups (see Theorems \ref{APgrowth} and \ref{Hyp}): it includes for example all 
  Coxeter groups \cite{Hum} and all Gromov hyperbolic groups \cite{Gro, GdH}. 
  All countable groups having the Haagerup property \cite{CCJJV}
  and having subexponential H-growth (w.r.t.\ some Haagerup function) are also contained in this class. 

On the other hand, the main result of \cite{Haa1} may be interpreted as saying that free groups 
have some Fej\'er-like summation property, involving only finitely supported multipliers. 
We conclude this paper by giving some sufficient conditions
for this Fej\'er property and its twisted versions to hold.  

\medskip 
The influence of Haagerup's seminal paper \cite{Haa1} on our work should be evident.
We have also benefitted from many of its follow-ups (like \cite{DCHa, CH, ER, Jol, JoVa, BF, BP, CCJJV, BrNi, Oz}). 
It should finally be noted that  Zeller-Meier's result for amenable groups mentioned above 
is valid in the more general setting of  twisted $C^*$-crossed products by discrete groups 
and that a  proof of  the analogue of Fej\'er summation 
for usual $C^*$-crossed products by an action of $\Relativi$ is given in \cite{Dav}. 
(For more on this, see \cite{BC2}).

\section{Preliminaries}
Throughout this article
$G$ denotes a discrete group and $e$ its identity element. 

\subsection{On twisted group operator algebras}
The basic reference on this subject is \cite{ZM} (see also \cite{PaRa, PaRa1, PaRa2}). We give here a short review. 
We follow standard terminology and notation in operator algebras, 
as may be found for example in \cite{Di, Di2, Dav}.

\begin{definition}
 A  (normalized) \emph{2-cocycle}  on $G$ with values in  $\Toro$ is a map $\sigma: G \times G \to \Toro$
such that
$$\sigma(g, h)\sigma(gh, k)=\sigma(h, k)\sigma(g, hk) \quad (g, h, k \in G)$$
$$\sigma(g, e)=\sigma(e, g)=1 \quad\quad (g \in G).$$
The set of all normalized 2-cocycles will be denoted by $Z^2(G,\Toro)$.
\end{definition}
Examples of 2-cocycles on $\Relativi^N$ were given in the Introduction. Up to a natural equivalence (irrelevant for our purposes) they
 are always of this form (\cite{Bac, BB}). For other examples, see  e.g.\  \cite{Klep} (for abelian groups), \cite{Pac} (for the integer Heisenberg group),  
\cite{How} (for Coxeter groups), \cite{MM1} (for Fuchsian groups).
 
Projective representations associated with 2-cocycles were first considered by I.\ Schur in the case of finite groups and by G.\ Mackey in the
general case (see e.g.\ \cite{Ma, Kl}). 
\begin{definition}  A $\sigma$-\emph{projective unitary representation} $U$ of 
$G$ on a (non-zero) Hilbert space $\mathcal{H}$ 
is a map from G into the group $\mathcal{U} (\mathcal{H})$ of
unitaries on $\mathcal{H}$ such that
$$ U(g)U(h) = \sigma(g,h) U(gh) \quad\quad (g, h \in G).$$
 \end{definition}
 We then have $U(e) = I_\mathcal{H}$ (the identity operator on $\mathcal{H}$) and 
$$U(g)^*=\overline{\sigma(g, g^{-1})}U(g^{-1}), \quad g \in G.$$

\begin{definition} Let  $\sigma \in Z^2(G,\Toro)$. The \emph{left 
regular $\sigma$-projective unitary representation} $\Lambda_\sigma$ of $G$ on $\ell^2(G)$ is defined  by
$$ (\Lambda_\sigma(g) \xi) (h) = \sigma(g, g^{-1}h) \xi(g^{-1} h), \ \ \xi \in
\ell^2(G),\ g, h \in G.$$ 
\end{definition} 

Choosing $\sigma$ to be the trivial 2-cocycle ($\sigma=1$)  gives  the  left 
regular representation of $G,$ which we will denote by $\lambda.$ 
Some authors prefer a unitarily equivalent definition of the left regular $\sigma$-projective
unitary representation  of $G$ (and others prefer right versions), but we have chosen to follow \cite{Pa}.

From now on, we fix $\sigma \in Z^2(G, \Toro)$. Letting $\{\delta_{h}\}_{h \in G}$ denote the canonical basis of $\ell^2(G),$
 we then have
 $$\Lambda_\sigma(g) \delta_{h} = \sigma(g,h) \delta_{gh}, \quad g, h \in G,$$
 so, especially,  $\Lambda_\sigma(g) \delta= \delta_{g}$, where $\delta= \delta_{e}$.
 
\begin{definition}
The \emph{reduced twisted group $C^*$-algebra} $C_{r}^{*}(G, \sigma)$ 
(resp.\ the \emph{twisted group von Neumann algebra} $vN(G, \sigma)$) is 
 the $C^*$-subalgebra (resp.\ von Neumann subalgebra) of $B(\ell^2(G))$ generated by the set
$\Lambda_{\sigma}(G),$ that is, the closure in the operator norm (resp.\ weak operator) topology of 
the $*$-algebra $\Complessi(G,\sigma) = $Span$(\Lambda_{\sigma}(G)).$  
\end{definition}

\begin{definition} We let $\tau$ denote the linear functional on $vN(G, \sigma)$ given by
$$\tau(x)= (x \delta, \delta), \; x \in vN(G, \sigma).$$
For $ x \in vN(G, \sigma)$, we set $\| x \|_{2}= \tau(x^* x)^{1/2}\, $ and $\ \widehat{x} = x \delta \in \ell^2(G)$.  
\end{definition} 

The following fundamental result is well known.
\begin{proposition}\label{LI}
The functional $\tau$ is a faithful, tracial state on $vN(G,\sigma)$  and
$\| \cdot \|_{2}$ is a norm on $vN(G,\sigma)$. 

Moreover, the map $x \rightarrow \widehat{x} $ is a linear isometry from 
 $(vN(G, \sigma), \| \cdot \|_{2})$ to $(\ell^2(G), \| \cdot \|_{2}),$ 
which sends  $\Lambda_{\sigma}(g)$ to $\delta_{g}$ for each $g \in G.$ 
\end{proposition}

\begin{definition} The value
 $\widehat{x}(g) \in \Complessi$ is called the \emph{Fourier coefficient} of $x \in vN(G,\sigma)$ at $g \in G.$ 
 \end{definition}
 To justify this definition, we first remark that $\tau$ corresponds to integration w.r.t.\ the normalized  Haar measure in the classical case.
Hence, we may  consider $\tau$ as the normalized "Haar functional" on $vN(G,\sigma)$. Then we have 
$$\widehat{x}(g) = (x \delta, \delta_{g})=(x \delta, \Lambda_{\sigma}(g) \delta) 
= \tau(x \Lambda_{\sigma}(g)^*) $$
for all $x \in vN(G,\sigma)$ and $g \in G$. Further, we also record that
$$ \| \widehat{x} \|_{\infty}   \leq  \| \widehat{x} \|_{2} = \| x \|_{2} \leq   \| x \|.$$

\begin{definition}The (formal) \emph{Fourier series} of $ x \in vN(G, \sigma)$ is the series 
$$\sum_{g\in G} \widehat{x}(g)\Lambda_{\sigma}(g).$$
\end{definition}

Note that this series does not necessarily converge in the weak operator topology (see \cite{Mer}). However, the
following result follows readily from Proposition \ref{LI}.
\begin{proposition}\label{L2} Let $ x \in vN(G, \sigma)$. Then 
$$  x = \sum_{g\in G} \widehat{x}(g)\Lambda_{\sigma}(g) \  \  
\; (\text{w.r.t.} \| \cdot \|_{2}).$$
\end{proposition}

\smallskip The Fourier series representation of
 $x \in vN(G, \sigma)$ is unique. 
 More generally, the following holds. 
 
 \begin{proposition} \label{FU} Let $\xi : G \to \Complessi$ and suppose that
the series $\sum_{g\in G}\xi(g) \Lambda_{\sigma}(g)$ converges to some $x \in vN(G, \sigma)$
w.r.t.\   $\| \cdot \|_{2}. $ Then  $\xi \in \ell^2(G)$ and $\xi = \widehat{x}.$
\end{proposition} 
\begin{proof} For any finite subset $F$ of $G,$  
set  $a_{F}=\sum_{g\in F}\xi(g) \Lambda_{\sigma}(g)$
and let $\chi_{F}$ denote the characteristic function of $F.$
Then we have $\widehat{a_{F}}= \xi \chi_{F} =: \xi_{F}.$ Now the assumption says that
$a_{F} \to x$ w.r.t.\  $\| \cdot \|_{2},$ which implies that $\xi_{F}  \to \widehat{x}$ in 
$\ell^2$-norm. This implies that $\sum_{g \in F}| \xi(g)|^2 \leq \| \widehat{x} \|_{2}^2$ for all finite subset $F$ of $G,$
hence $\xi \in \ell^2(G).$ But then $\xi_{F} \to \xi$ in $\ell^2$-norm and we get $\widehat{x}= \xi.$ 
\end{proof}

\begin{definition} We set 
$$CF(G,\sigma):= \{ x \in C^*_{r}(G,\sigma) \, | \, 
\sum_{g\in G} \widehat{x}(g)\Lambda_{\sigma}(g) \
\text{ is convergent
in operator norm} \}.$$
\end{definition}
\begin{proposition}If  $x \in CF(G,\sigma)$, then its Fourier series necessarily converges to $x$ in operator norm.
\end{proposition}
\begin{proof} This follows from Proposition \ref{L2}
and the fact that $\| \cdot \|_{2} \leq \| \cdot \|$.
\end{proof}

Let $f \in \ell^1(G).$ The series $\sum_{g\in G} f(g) \Lambda_{\sigma}(g)$ is clearly 
absolutely convergent in operator norm and we shall denote its sum by $\pi_{\sigma}(f).$ Then we have
$\| \pi_{\sigma}(f) \| \leq  \| f \|_{1}$ and 
$$\widehat{\pi_{\sigma}(f)} = (\sum_{g\in G} f(g)\Lambda_{\sigma}(g)) \delta = 
\sum_{g\in G} f(g)\delta_{g} = f.$$ 
Note that in the sequel, we will use the more suggestive notation $\pi_{\lambda}$ instead of $\pi_{1}$ (since we write $\lambda$ 
instead of $\Lambda_{1}$).

Let now $x \in vN(G,\sigma)$ and assume that $\widehat{x} \in \ell^1(G).$ Then we get
$ \widehat{\pi_{\sigma}(\widehat{x})}= \widehat{x},$ hence $\pi_{\sigma}(\widehat{x}) = x. $
Therefore, $x \in CF(G,\sigma)$ and $ \| x \| = \| \pi_{\sigma}(\widehat{x})Ê\| \leq \| \widehat{x} \|_{1}$.

Summarizing, we get the following.
\begin{proposition}
$$\{ x \in vN(G,\sigma) \; | \; \widehat{x} \in \ell^1(G) \} = \pi_{\sigma}(\ell^1(G)) \subseteq CF(G,\sigma).$$
\end{proposition}

\smallskip Twisted group operator algebras may alternatively be described with the help  of twisted convolution. 
\begin{definition} Let $\xi, \eta \in \ell^2(G).$ The complex function $\xi \as \eta$ on $G$ given by
 $$(\xi \as \eta)(h) = \sum_{g\in G} \xi(g) \sigma(g, g^{-1}h) \eta(g^{-1}h), \ h \in G$$
is called the  $\sigma$-\emph{convolution product}  of $\xi$ and $\eta$.
\end{definition}
As $| (\xi \as \eta)(h) | \leq (|\xi| \ast |\eta| ) (h), \ h \in G,$
it is straightforward to verify that $\xi \as \eta$ is a well defined bounded function on $G$ satisfying
$$\| \xi \as \eta \|_{\infty} \leq \| \,  |\xi| \ast |\eta| \, \|_{\infty} \leq \| \xi \|_{2} \| \eta \|_{2}.$$
Notice also that  $\delta_{a} \as \delta_{b} = \sigma(a,b) \delta_{ab}, \, a, b \in G.$ 

\vspace{1ex}
\noindent One can now check that if $x\in vN(G, \sigma)$ and $\eta \in \ell^2(G)$, then $x \eta = \widehat{x} \as \eta.$
 The usual properties of convolution carries over to twisted convolution. For example, we have 
\begin{proposition} Let $p \in \{Ê1, 2\}$, $f \in \ell^1(G), \eta \in \ell^p(G).$ Then $ f \as \eta \in \ell^p(G)$ and
$$\| f \as \eta \|_{p} \leq \| f \|_{1} \|\eta\|_{p}$$
 Moreover, the Banach space $\ell^1(G)$ is a
  Banach $*$-algebra,   denoted by $\ell^1(G,\sigma)$,  with respect to twisted convolution and involution  given by 
  $$f^*(g)= \overline{\sigma(g,g^{-1})}\, \overline{f(g^{-1})}, \, g\in G.$$ 
\end{proposition}

As $\pi_{\sigma}(f) \eta = f \as \eta$ whenever $ f \in \ell^1(G) $ and $ \eta \in \ell^2(G), $
the map $\pi_{\sigma} : \ell^1(G) \to C^*_{r}(G,\sigma)$ is easily seen to be a faithful $*$-representation of $\ell^1(G, \sigma)$ 
 on $ \ell^2(G)$, and $C_{r}^*(G, \sigma)$ is the closure of $\pi_{\sigma} (\ell^1(G))$ 
 in the operator norm.
 Moreover, there is a bijective correspondence $U \rightarrow \pi_{U}$ between $\sigma$-projective 
 unitary representations of $G$ and non-degenerate $*$-representations of $\ell^1(G,\sigma)$ determined by
 $$\pi_{U}(f)= \sum_{g \in G}f(g)U(g), \ f \in \ell^1(G),$$
 (the series above being obviously absolutely convergent in operator norm), the inverse correspondence 
 being given by $U_{\pi}(g)= \pi(\delta_{g}), g \in G.$ One may then pass to the
the enveloping C$^*$-algebra \cite{Di} of  $\ell^1(G, \sigma)$, which is denoted by $C^*(G, \sigma)$ and
 called the \emph{full twisted group}  $C^*$-\emph{algebra associated to} $(G,\sigma).$ When $G$ is amenable, 
the extension of $\pi_{\sigma}$ to $C^*(G, \sigma)$ is faithful (\cite{ZM}), and $C^*(G, \sigma)$ may then be identified
with $C_{r}^*(G, \sigma)$ via this isomorphism. 
 
\subsection{On amenability, Haagerup property and length functions}

\begin{definition} The group $G$ is called \emph{amenable} if there 
exists a left translation invariant state on $\ell^{\infty}(G).$ 
\end{definition} 
Amenability of $G$ can be formulated in a huge number of
equivalent ways (see \cite{Di, Pat, Pi, Wa}). We will make use of the following characterizations. 
As usual, a complex function $\varphi$ on $G$ is called \emph{normalized} when $\varphi(e)=1$.

\begin{theorem} \label{A}
The group $G$ is amenable if and only if one of the following conditions holds :
\begin{flushleft}
1) $G$ has a \emph{F\o lner net} $\{F_\alpha\}$, that is,  
 each $F_\alpha$ is a finite non-empty subset of $G$ and we have
\begin{equation}
\label{folner}
\frac{|g F_\alpha \mbox{\scriptsize${\triangle}$}  F_\alpha|}{|F_\alpha|} 
\underset\alpha\to 0 
\quad \text{ for every}\ g \in G . 
\end{equation}

\smallskip 2) There exists a net $\{ \varphi_{\alpha}Ê\}$ of normalized positive definite functions on $G$ with finite support such that
 $\varphi_{\alpha} \to 1$ pointwise on $G.$ 
  
\smallskip 3) There exists a net $\{ \psi_{\alpha}Ê\}$ of normalized positive definite functions in $\ell^{2}(G)$  
such that $\psi_{\alpha} \to 1$ pointwise on $G.$ 

\smallskip 4) $\Big| \, \sum_{g \in G}f(g) \,\Big|  \leq \| \sum_{g \in G} f(g) \lambda(g)\|$ \; 
for all $f \in \ell^1(G).$

\end{flushleft}
\end{theorem}
In the sequel, we will write p.d.\ instead of positive definite. In the same way, we will write n.d.\ instead of negative definite
(we follow here \cite{BCR}; n.d.\ functions are called conditionally negative definite by some authors).

\medskip A weakening of  2) in Theorem \ref{A} leads to the following concept (see \cite{CCJJV}). 
\begin{definition} The group $G$ is said to have the \emph{Haagerup property} (or to be a-T-menable) if there exists a net 
$\{ \varphi_{\alpha}\}$ of normalized p.d.\ functions vanishing at infinity on $G$ (that is, 
$\varphi_{\alpha} \in c_{0}(G)$ for all $\alpha$) and converging pointwise to $1.$ 
\end{definition}

Clearly, all amenable groups have the Haagerup property. All free groups also have this property, as first established in \cite{Haa1}. 
We refer to \cite{CCJJV} for other examples, as well as many characterizations of the Haagerup property. We will need the following one.

\begin{proposition} \label{H} Assume that $G$ is countable. 
Then $G$ has the Haagerup property if and only if  there exists a \emph{n.d}.\ function 
$h: G \to [0, \infty)$ which is \emph{proper}, that
is, \ $h^{-1}([0, t])$ is finite for all $t \geq 0$.
We will call such a function $h$ a \emph{Haagerup function} on $G.$ 
\end{proposition}

Concerning n.d.\ functions, we recall the following result 
of Schoenberg (which is used to prove
 Proposition \ref{H}; see \cite[Theorem 2.2]{BCR} for a more general statement). 

 \begin{theorem} \label{S}
 A function $\psi : G \to \Complessi$ is
n.d.\ if and only if $e^{-t \psi}$ is p.d. for all $t >0$ (equivalently, $r^\psi$ is p.d for all $0 < r < 1$). 
\end{theorem}

An interesting class of functions on a group is the class of length functions (see e.g.\ \cite{Con2, Jol, JoVa}). 
\begin{definition} A function $L: G \to [0, \infty)$ is called a \emph{length} function if 
$L(e)= 0, L(g^{-1}) = L(g)$ and $L(gh)\leq L(g)+L(h)$ for all $g, h \in G.$
\end{definition}
If $G$ acts isometrically on a metric space $(X,d)$ and $x_{0} \in X,$ then $L(g):=d(g\cdot x_{0}, x_{0})$ 
gives a ``geometric'' length function on $G.$ (All length fuctions can be described in this way). If $G$
is finitely generated and $S$ is a finite generator set for $G$, then the word-length function $g \to L_{S}(g)$ 
(w.r.t.\   to the letters from $S \cup S^{-1}$) gives a length function on $G$, which we will call \emph{algebraic}. 
All such algebraic length functions are equivalent in a natural way. 

\bigskip
Length functions may be used to define growth conditions. The reader should consult  \cite{Har, Jol, Pat, Wa} for more details.
\begin{definition}
Let $L$ be a length function on $G$. For $r \in \Reali^{+}$, set \\
 $B_{r,L}:=\{\, g\in G \; | \; L(g) \leq r  \, \}.$ Then one says that  

 \smallskip 
 $G$ has
\emph{polynomial growth (w.r.t.\  $L$)}  if there exist  $K, p > 0 $ such that $|B_{r, L}| \leq K (1 + r)^p$ for all $r\in \Reali^{+}$, 

\smallskip  $G$ has \emph{subexponential growth (w.r.t.\  $L$)} 
if for any $b >1,$ there is some $r_{0} \in \Reali^{+}$ such that 
 $|B_{r,L}| < b^r$ for all $r \geq r_{0}.$ 
 \end{definition}
 
\begin{definition} Let $G$ be finitely generated.  Then $G$ has \emph{polynomial} (resp.\  \emph{subexponential}) \emph{growth} 
if it has polynomial (resp.\ subexponential) growth w.r.t.\  
some (or, equivalently, any) algebraic length on $G.$ 
\end{definition}
Note that if $G$ is finitely generated and has polynomial (resp.\ subexponential) growth w.r.t.\  to some  length function $L$ on $G,$ then 
$G$ has polynomial (resp.\ subexponential) growth. 
In addition, we mention:

\begin{theorem} \label{GG} Let $G$ be finitely generated and let $S$ be a generator set.
\begin{flushleft}
1) If $G$ has polynomial growth, then $\{B_{k, L_{S}}\}_{k \geq 0}$ is a F\o lner sequence for $G$ (see \cite{Har}).

\smallskip 2) If $G$  has subexponential growth, then there is a subsequence of $\{B_{k, L_{S}}\}_{k \geq 0}$ which is a F\o lner sequence for $G$
(see \cite{Har}).

\smallskip  3)  $G$ has polynomial growth if and only if $G$
is almost nilpotent (see \cite{Pat, Wa}). 

\smallskip 4)  $G$ can have subexponential growth without having polynomial growth (see \cite{Pat, Wa}). 
\end{flushleft}
\end{theorem}

 Length functions arise naturally in connection with the Haagerup property.
 \begin{proposition}\label{HL}
Let $G$ be countable. Then $G$ has Haagerup property if and only if it has a Haagerup length function. 
\end{proposition}
 \begin{proof} Assume that $G$ has the Haagerup property.
Let $h$ be a Haagerup function on $G$ (cf.\ Theorem \ref{H}). Then, as $h$ is n.d.,  $L=h^{1/2}$ is also n.d.\ (see \cite[Corollary 2.10]{BCR}).
 Further, $L$ is clearly proper. Finally, $L$ is a length function on $G$ : 
this follows from \cite[Proposition 3.3]{BCR} (the standing assumption that $G$ is abelian is not
 used in the proof of this proposition). Hence $L$ is a Haagerup length function on $G.$ The converse implication is trivial.
 \end{proof}

 In some cases,  Haagerup length functions are geometrically given :
 this happens for example when $G$ acts isometrically and metrically properly on a tree, or on a $\Reali$-tree,   
\cite{Bo, Va}. In the case of finitely generated groups, Haagerup length functions are often algebraically given : 
 this is at least true for free abelian groups \cite{BCR}, free  groups \cite{Haa1, CCJJV} and 
 Coxeter groups \cite{BJS}. 

\section{Convergence of Fourier series and decay properties}
Throughout the rest of this paper, we let $\sigma \in  Z^2(G, \Toro)$ and denote by $\K(G)$ the set of all complex
functions on $G$ having finite support. 

\begin{definition}
Let $\L$ be a subspace of $\ell^2(G)$ which contains $\K(G),$ 
let $\| \cdot \|'$ be a norm on $\L $ and   $ \xi \in \L.$ When
$F$ is finite subset of $G,$ set $\xi_{F} = \xi \chi_{F},$ where
$\chi_{F}$ denotes the characteristic function of $F.$ 

We say that $\xi \to 0$ at infinity w.r.t.\  $\| \cdot \|'$ if 
for every $\varepsilon > 0,$ there exists a finite subset $F_{0}$
of $G$ such that $\| \xi_{F} \|' < \varepsilon$ for all
finite subsets $F$ of $G$ which are disjoint from $F_{0}.$

\end{definition}

\begin{definition}\label{LDP}
Let $\L$ be a subspace of $\ell^2(G)$ which contains $\K(G).$ We say that
$(G, \sigma)$ has the $\L$-decay property (w.r.t.\  $\| \cdot \|'$) if there exists a norm $\| \cdot \|'$ on $\L$
such that the following two conditions hold:
\begin{itemize}
\item[i)] For each $\xi \in \L$ we have $\xi \to 0$ at infinity w.r.t.\  $\| \cdot \|'.$

\item[ii)] The map $f \to \pi_{\sigma}(f)$ from $(\K(G), \| \cdot \|' )$ to $(C_{r}^*(G, \sigma), \| \cdot \|)$
is bounded.
\end{itemize}
We will simply say that $G$ has the $\L$-decay property (w.r.t.\  $\| \cdot \|'$)  if $(G, 1)$ has 
the $\L$-decay property (w.r.t.\  $\| \cdot \|'$).
\end{definition}

Due to the following proposition, it is sufficient to establish decay properties only for $G$ in all natural cases
we are aware of.

\begin{proposition} \label{norm} Assume that $G$ has the $\L$-decay property w.r.t.\   $\| \cdot \|'$ and that
$\| \, |f | \, \|'  = \| f \|'$ for all $f \in \K(G).$ Then $(G, \sigma)$ has the $\L$-decay property w.r.t.\   $\| \cdot \|'$.
\end{proposition}
\begin{proof} 
Let  $C > 0$ be the norm of the map
$f \to \pi_{\lambda}(f)$ from $(\K(G), \| \cdot \|' )$ to $(C_{r}^*(G), \| \cdot \|)$. Let $f \in \K(G)$ and $\eta \in \ell^2(G).$ Then
\begin{align*}
 \| \pi_{\sigma}(f) \eta \|_{2} & 
= \| f \as \eta \|_{2} \leq  
 \| \,  | f |  \ast  | \eta  | \,  \|_{2} 
= \|\pi_{\lambda}(| f |) | \eta | \, \|_{2} \\
& \leq C \| \, | f | \, \|'  \, \| \, | \eta | \, \|_2
 = C \| f \|' \, \|\eta \|_{2}
\end{align*}
 Hence, we have 
$ \|\pi_{\sigma}(f) \| \leq  C \| f \|'$ for all $f \in \K(G).$ As the first condition in Definition \ref{LDP}
is independent of $\sigma,$ the assertion follows.
\end{proof}
The above proposition has previously been established by I. Chatterji \cite{Cha} (in a special situation).

\begin{lemma}\label{CFS1}
Assume that $(G, \sigma)$ has the $\L$-decay property w.r.t.\  $\| \cdot \|'$. 

Let $\xi \in \L.$
Then the series $\sum_{g \in G} \xi(g) \Lambda_{\sigma}(g)$ converges in operator norm to some 
$a \in C_{r}^*(G, \sigma)$ satisfying $\widehat{a} = \xi.$ We will denote this $a$ by $\tilde{\pi}_{\sigma}(\xi).$

Letting $\tilde{\pi}_{\sigma} : \L \to C_{r}^*(G, \sigma)$ denote the associated map, we then have
$\tilde{\pi}_{\sigma} (\L) \subseteq CF(G, \sigma)$.
\end{lemma}
\begin{proof}Using that ii) holds, we get that there exists  $C >0$ such that
$$\| \sum_{g\in F} \xi(g) \Lambda_{\sigma}(g) \| = \| \pi_{\sigma}(\xi_{F}) \| 
\leq  C \| \xi_{F}\|' $$ for any finite subset $F$ of $G.$ Now, using that i) holds, we deduce then  immediately that the net 
$\{ \sum_{g\in F} \xi(g) \Lambda_{\sigma}(g) \}_{F}$, indexed over the finite subsets of $G$ ordered by inclusion,
satisfies the Cauchy criterion \cite[9.1.6]{Di3} w.r.t.\  operator norm. Hence this net converges in operator norm to some $a \in C_{r}^*(G, \sigma).$
But then it also converges to $a$ w.r.t.\  $\| \cdot \|_{2},$  hence we have $ \widehat{a} = \xi$ by Proposition \ref{FU}, as desired. 
The last statement follows immediately.
\end{proof}

\begin{theorem}\label{CFS2}
Assume that $(G, \sigma)$ has the $\L$-decay property w.r.t.\  $\| \cdot \|'$. \\
Set $ \L^{\vee}(G, \sigma) = \{ x \in vN(G, \sigma)  \: |  \: \widehat{x} \in \L \} .$
Then $$\L^{\vee}(G, \sigma) = \tilde{\pi}_{\sigma} (\L) \subseteq CF(G, \sigma).$$
\end{theorem}
\begin{proof}Let $ x  \in \L^{\vee}(G, \sigma).$ From  Lemma \ref{CFS1} (with $\xi = \widehat{x}$),
 we get that the Fourier series of $x$
converges in operator norm to $ \tilde{\pi}_{\sigma} (\widehat{x}) \in C_{r}^*(G, \sigma)$ and that 
 $\widehat{\tilde{\pi}_{\sigma} (\widehat{x})} = \widehat{x}.$ By uniqueness, this implies that 
$\tilde{\pi}_{\sigma} (\widehat{x}) = x.$ Thus we have shown that $\L^{\vee}(G, \sigma) \subseteq \tilde{\pi}_{\sigma} (\L)$.
 We also know that $\tilde{\pi}_{\sigma} (\L) \subseteq CF(G, \sigma)$ from Lemma \ref{CFS1}.

Finally, if $x \in \tilde{\pi}_{\sigma}(\L),$ so $x = \tilde{\pi}_{\sigma}(\xi)$ 
for some $\xi \in \L,$ then  Lemma \ref{CFS1} says that $\widehat{x} = \xi \in \L.$ Hence, 
$\tilde{\pi}_{\sigma} (\L) \subseteq \L^{\vee}(G, \sigma) $. 
\end{proof}

It is almost immediate that $(G, \sigma)$ has the $\ell^1(G)$-decay property w.r.t.\  $\| \cdot \|_{1}.$ Anyhow, 
we already saw in Section 2 that the assertions in Lemma \ref{CFS1} and Theorem \ref{CFS2} hold when $\L = \ell^1(G).$

\vspace{1ex}
As another source of examples, we shall now consider weighted spaces. We 
establish first some notation.

\vspace{1ex}
Let $\kappa: G \to [1,\infty), 1 \leq p \leq \infty$ 
and define $$\L_\kappa^p 
= \{\xi : G \to \Complessi \ | \ \xi\kappa \in \ell^p (G)\}
\subseteq \ell^p (G),$$ which becomes a Banach space w.r.t.\  the
norm  $\|\xi\|_{p,\kappa} = \|\xi \kappa\|_p$.
Clearly, $\L_\kappa^p \subseteq \L_\kappa^q $ and $\| \cdot \|_{q,\kappa} \leq \| \cdot \|_{p,\kappa}$ whenever
 $1 \leq p \leq q \leq \infty, $ 
 while $\L_\gamma^p \subseteq \L_\kappa^p $ and $\| \cdot \|_{p,\kappa} \leq \| \cdot \|_{p,\gamma}$ whenever
 $\gamma: G \to [1,\infty)$ is such that $\kappa \leq \gamma .$

\begin{definition}
We say that $(G,\sigma)$ (resp.\ $G$)
is $\kappa$-decaying if $(G,\sigma)$ (resp.\ $G$) has the $\L_{\kappa}^2$-decay property w.r.t.\  $\|\cdot\|_{2,\kappa}.$
\end{definition}

\begin{proposition}  \label{L2D}    

\begin{flushleft} 
\noindent 1) $G$ is $\kappa$-decaying if and only if 
the linear map $f \to \pi_{\lambda}(f)$ from $(\K(G), \| \cdot \|_{2,\kappa} )$ to $(C_{r}^*(G), \| \cdot \|)$ is bounded. 
The norm of this map will then be called the $\kappa$-\emph{decay constant} of  $G$. 

\smallskip 
\noindent 2) Assume that $G$ is $\kappa$-decaying. Then $(G, \sigma)$ is $\kappa$-decaying and we have
$$ \{ x \in vN(G, \sigma) \: | \: \widehat{x} \in \L_{\kappa}^2 \} = \tilde{\pi}_{\sigma} (\L_{\kappa}^2) \subseteq CF(G, \sigma).$$
\end{flushleft}
\end{proposition}
\begin{proof}

\noindent 1) The fact that $\L_{\kappa}^2$ satisfies condition i) in Definition \ref{LDP} 
w.r.t.\  $\|\cdot\|_{2,\kappa}$ is elementary classical analysis. 

\noindent 2) As  $\| \, |f| \, \|_{2,\kappa}= \| \, f\|_{2,\kappa}$ for all $f \in \K(G)$, the first assertion follows from Proposition  \ref{norm}. The second
is then a consequence of Theorem \ref{CFS2}.
\end{proof}

\begin{example} \label{IS} Assume that $G$ is countable and  $\kappa$ satisfies condition
(IS), by which we mean that $\kappa^{-1} \in \ell^2(G).$ 

\smallskip \noindent  Then the Cauchy-Schwarz inequality immediately
gives that $\L^2_{\kappa} \subseteq \ell^1(G) $ and
 $\| f \|_{1} \leq \| \kappa^{-1} \|_{2} \, \| f \|_{2, \kappa}, f \in \L^2_{\kappa}$. 
  As $\| \pi_{\lambda}(f) \| \leq  \| f \|_{1}$ for all  $f \in \K(G),$  we get 
 $$\| \pi_{\lambda}(f) \| \leq \| \kappa^{-1} \|_{2} \, \| f \|_{2, \kappa}, \; \; f \in  \K(G).$$
 Hence,  $G$ is $\kappa$-decaying (with decay constant at most  $\| \kappa^{-1} \|_{2}$). However, note that in such a case,
 the conclusion of Proposition \ref{L2D}, part 2), brings nothing new as $\pi_{\sigma}(\L^2_{\kappa}) \subseteq 
 \pi_{\sigma}(\ell^1(G)) $.
 \end{example} 
  
 More concretely, assume that $G$ is finitely generated and let $L$ denote any algebraic length function on $G.$ 
For $t> 0, $ set $\kappa_{t}= \exp(tL^2).$ Then $\kappa_{t}$ satisfies  (IS) (see e.g.\ the proof of   
\cite[Proposition 24]{Con2}).
One may also consider  $\gamma_{a}=a^{L}, a >1.$
Then  $\gamma_{a}$ is easily seen to satisfy (IS) for all $a>1$ 
whenever $G$ has subexponential growth. Hence, $G$ is $\gamma_{a}$-decaying for all $a >1$ in this case. 
 As we will see later in this section, the same conclusion can still be drawn for many nonamenable groups. 

\bigskip The case 
$\kappa_{L, s} = (1 + L)^s $ where $L$ is a length function  and 
$s > 0$   has received a lot of attention in connection with the rapid decay property for groups, introduced in \cite{Jol}.
Using our terminology, $G$ has the RD-property (w.r.t.\  a length  function $L$)  if and only if there exists some  $s_{0} > 0$ such
that $G$ is $\kappa_{L, s_{0}}$-decaying. Note that when $G $ is amenable, $G$ has the RD-property (w.r.t.\ $L$) 
if and only if $G$ has polynomial growth (w.r.t.\   $L$) (see \cite[Corollary 3.1.8]{Jol} and \cite{V, CR}).
When $G$ is finitely generated, one just talks about the RD-property, having in mind 
that $L$ is then chosen to be any algebraic length function on $G.$  

 Much of the interest around the RD-property is due to the following : 
when $G$ has the RD-property (w.r.t.\ $L$), then the canonical image of the Fr\'echet space 
 $H_{L}^{\infty} = \cap_{s > 0}\,  \L^2_{\kappa_{L, s}}$ (w.r.t.\   the obvious family of seminorms),
which is thought as representing a space of "smooth" functions on the "dual" of $G,$ is a dense spectral (= inverse-closed) $*$-subalgebra of $C^*_{r}(G).$ 
For more about this and the RD-property, see  e.g.\ \cite{Jol2, JoVa, JS, V, Cha, CR, Lu2} and references therein. See also the end of this section.
 
\bigskip

Let now $L$ denote the usual word-length function on a free group $\mathbb{F}_{n}.$ It follows from \cite{Haa1} 
(see  also \cite{Jol,V}) that $\mathbb{F}_{n}$ is $\kappa_{L,2}$-decaying, hence $\mathbb{F}_{n}$ has the RD-property. 
This may be seen as a consequence of the fact that $\mathbb{F}_{n}$ has "polynomial H-growth".  
To explain this, we begin with a fundamental lemma.

\begin{lemma}\label{Clem}

Let  $E$ be a  non-empty finite subset of $G.$ Define
$$c(E) = 
\sup \{ \, \|\pi_{\lambda}(f)\| \ | \  f \in \K(G), \, {\rm supp}(f) \subseteq E,  
\|f\|_2 = 1\}.$$
Then $1 \leq c(E) \leq |E|^{1/2}$.

\medskip
\noindent If $G$ is amenable, then $c(E) = |E|^{1/2}.$
\end{lemma}

\begin{proof}
Note first that if  $a \in E,$ then $ \|\delta_{a} \|_{2} = 1$ and  
$\|\pi_{\lambda}(\delta_{a})\| = \|\lambda(a)\|=1.$ Hence
 $c(E) \geq 1.$ 

Next, we have
$$\|\pi_1(f)\| \leq \|f\|_1 = \sum_{g \in E} |f(g)|
\leq |E|^{1/2} (\sum_{g \in E} |f(g)|^2)^{1/2} = |E|^{1/2} \|f\|_2$$
for every $f \in \K(G)$ with ${\rm supp}(f) \subseteq E$.
So $c(E) \leq |E|^{1/2}.$

Finally, assume that $G$ is amenable. Set $f = (1/ |E|^{1/2}) \chi_{E}.$ Then we have $\|f\|_2 = 1$ and
$ |E|^{1/2} = \| fÊ\|_{1} = \|\pi_{\lambda}(f)\|$ (cf.\ Theorem \ref{A}, part 4). Hence we get $ |E|^{1/2} \leq  c(E)$ 
and the last assertion follows.
\end{proof}
Obviously, $c(E)$ is what we called the \emph{Haagerup content} of $E$ in the Introduction. 
We leave to the reader to check that $c(E) \leq c(F) $ whenever $E \subseteq F$ and 
$c(E \cup F) \leq c(E) + c(F)$ whenever $E$ and $F$ are pairwise disjoint ($E, F$ being finite nonempty subsets of $G$).

The computation of $c(E),$ or just finding an upper bound for it  
better than $|E|^{1/2}$ when $G$ is nonamenable, appears to be quite challenging in general. 
It has been dealt with in some special cases (e.g.\ \cite{L, AO, Flo, Haa1, FTP, Coh, JoVa, Fen}), often in connection 
with the related problem of estimating the norm $\|\pi_{\lambda}(f)\|$ for $f \in \K(G)$ (especially when $f = \chi_{E}$). 

We can now measure "H-growth" instead of growth by using the Haagerup content instead of 
the square root of cardinality for finite subsets of $G$. 

\begin{definition}\label{pHgrowth}
Let $G$ be countable and $L : G \to [0, \infty)$ be proper. Set $B_{r, L}:=\{ g \in G  \, |  \,Ê L(g) \leq r \}$ for each $r \in \Reali^{+}$.
Then we say that 

\smallskip $G$ has \emph{polynomial H-growth (w.r.t.\  $L$)} if there exist  $K, p > 0$ such that 
$c(B_{r, L}) \leq K (1+ r)^p$ for all $r \in \Reali^{+}$. 

\medskip
 \noindent Further, we say that $G$ has \emph{subexponential H-growth (w.r.t.\   $L$)}
 if for any $b >1,$ there exists some $r_{0} \in \Reali^{+}$ such that 
$c(B_{r, L}) < b^r $ whenever $r\geq r_{0}.$

\end{definition}
It is clear from Lemma \ref{Clem} that when $G$ is amenable and $L$ is a proper length function on $G$, then
polynomial (resp.\ subexponential) H-growth (w.r.t\ $L$) reduces to polynomial (resp.\ subexponential) growth (w.r.t\ $L$). 

Using the properties of the Haagerup content mentioned above, one cheks without trouble the following useful lemma. 

\begin{lemma}\label{PSE} Let $G$ be countable and $L : G \to [0, \infty)$ be proper.
\smallskip \noindent For each  $ k \in \Relativi,$ $k \geq 0$, set $A_{k, L}= \{ g \in G  \, |  \,Êk \leq L(g) < k + 1 \}$ 
  and 
  
\smallskip \noindent $C_{L}(k)= c(A_{k, L}) \  \text{if} \ A_{k, L} 
\ \text{is nonempty,  $C_{L}(k)=0$ otherwise}.$

\smallskip Then $G$ has polynomial H-growth if and only if there exist constants $K, p > 0$ such 
 that $C_{L}(k) \leq K (1+ k)^{p}$ for all $k \geq 0.$

\smallskip Further, $G$ has subexponential H-growth if and only if for any $b>1,$ there exists some 
$k_{0} \in \Naturali $ such that $C_{L}(k) < b^k$ whenever $k \geq k_{0}$.
\end{lemma}  

\begin{example}
\label{Hex}
Using Lemma \ref {PSE}, a careful look into the existing literature provides us with many examples
of (nonamenable) groups having polynomial H-growth. 
\begin{itemize}
\item[1)]
Let $G = {\mathbb F}_n, n < \infty,$  denote a free group and let $L$ denote the natural algebraic length on $G.$
Then we have $C_{L}(k) \leq k+1$ for all $k \geq 0$ (see \cite{Haa1} for $n=2$ and \cite{V}
for a nice geometric proof of the general case due to T. Steger). Hence $G$ has polynomial H-growth (w.r.t.\   $L$).
\item[2)] More generally, let $G$ denote a  Gromov hyperbolic group \cite{Gro, GdH} 
and let $L$ denote some algebraic length on $G$. Then $G$ has polynomial H-growth (w.r.t.\   $L$). 
This may be deduced from \cite{Jol, dH2} (see also \cite{OzRi} and \cite{Co}) : in the course of the proof 
that $G$ has the RD property, it is implicitely shown that there exists a constant $K>0$ such that 
$C_{L}(k) \leq K (1+k)$ for all $k\geq 0.$ 
\item[3)]
Let $(G, S)$  denote a Coxeter group \cite{Hum} and let $L$ denote the word-length on $G$ 
(w.r.t.\   $S$). Then $C_{L}(k) \leq K(1+k)^{\frac{3}{2}P}$ for some $K > 0$ and 
$P \in \Naturali$, see \cite{Fen}. Hence $G$ has polynomial H-growth (w.r.t. $L$). 
Note that $G$ is nonamenable whenever it is neither finite nor affine \cite{dH}.
\item[4)] Let $G = G_{1}\ast_{A} G_{2}$ be an amalgamated free product  of groups  
 and let $L$ denote the "block" length on $G$ induced
by some integer-valued length functions $L_{j}$ on $G_{j}$ satisfying $L_{j}= 0$ on $A, \, j=1,2,$ (cf.\  \cite{Pic, Bo}).  
If $A$ is finite and each $G_{j}$ has polynomial $H$-growth (w.r.t.\   $L_{j}$), $j=1,2,$ then, adapting the proof 
of \cite[Theorem 2.2.2 (1)]{Jol}, one can deduce that $G$ has polynomial
H-growth (w.r.t.\   $L$). 
\end{itemize}
\end{example}

\medskip
To produce an example of  a nonamenable group $G$ which has subexponential, but not polynomial, H-growth 
(w.r.t.\ a length function $L$), one may proceed as follows. Pick  any finitely generated group $\Gamma$ which has subexponential, 
but not polynomial growth, and let $L_{1}$ denote some word-length function on $\Gamma.$ Then set  $G:= \Gamma \times \mathbb{F}_{2}$
and let $L$ be defined on $G$ by $L(g_{1},g_{2})= L_{1}(g_{1})+L_{2}(g_{2}),$  $L_{2}$ denoting the usual word-length function
on $\mathbb{F}_{2}.$ 

\medskip
Our interest in H-growth lies in the following.

\begin{theorem}

\label{Hgrowth}
Let $G$ be countably infinite and  $L : G \to [0, \infty)$ be  proper.  

\medskip \noindent
1) Assume that $G$  has polynomial H-growth (w.r.t.\ $L$).
Then there exist some $s_{0}>0$ such that $(G, \sigma)$ is $(1+L)^{s_{0}}$-decaying. 
Especially, if $L$ is a length function, then $G$ has the $\sigma$-twisted RD-property (w.r.t.\   $L$).

\smallskip \noindent
2) Assume that $G$ has subexponential H-growth (w.r.t.\ $L$).
Then $(G, \sigma)$ is $a^{L}$-decaying for all $a > 1.$ 
\end{theorem}

To prove this theorem, we will use the following.

\begin{lemma}
\label{partition}
Assume that $G$ is countably infinite and let $\{ E_{j}\}^\infty_{jÊ=0}$ be a partition of $G$ into
finite nonempty subsets.  

\smallskip \noindent Set $c_{j} = c(E_j), \,  j \geq 0$. Pick  $ d_{j} \geq 1$ for each  $j\geq 0$ such 
that $\sum_{j=0}^\infty (\frac{c_{j}}{d_{j}})^2 < \infty$.

\smallskip \noindent Define $\kappa:G \to [1, \infty)$ by $\kappa= \sum_{j=0}^{\infty} d_{j}\chi_{E_{j}}.$

\medskip
\noindent Then $G$ is $\kappa$-decaying. 
\end{lemma}

\begin{proof} 
Set $\chi_j = \chi_{E_j}, \ j \geq 0$. For $f \in \K(G)$, we have
\begin{align*}
\|\pi_{\lambda}(f)\| & = \|\sum_{j=0}^\infty \pi_{\lambda} (f \chi_j)\|  
\leq  \sum_{j=0}^\infty \|\pi_{\lambda} (f \chi_j)\| \\
& \leq \sum_{j=0}^\infty c_{j} \| f  \chi_j\|_2  
 =  \sum_{j=0}^\infty \frac{c_{j}}{d_j} d_j \| f  \chi_j\|_2  \\
& \leq \Big(\sum_{j=0}^\infty \Big(\frac{c_{ j}}{d_j} \Big)^2\Big)^{1/2}
\Big(\sum_{j=0}^\infty d_j^2 \| f  \chi_j\|_2^2 \Big)^{1/2}
= C \|f\|_{2,\kappa},
\end{align*}
where $C = (\sum_{j=0}^\infty (\frac{c_{j}}{d_{j}})^2)^{1/2}$. Hence $G$ is $\kappa$-decaying. 
\end{proof}
Note that if one also assumes that $G$ is amenable in this lemma, then one realizes easily that $\kappa$ satisfies (IS), so that
the assertion is  essentially trivial in this case (cf.\ Example \ref{IS}).

\medskip \emph{Proof of Theorem \ref{Hgrowth}}.
For each $k \geq 0,$ let $A_{k, L}$ and $C_{L}(k)$ be defined as in Lemma \ref{PSE}. 

\smallskip
Define $I = \{ k \in \Naturali \cup \{0 \}\, | \, A_{k, L} \ \text{is nonempty} \}$ and let 
$\{ k_{j}Ê\}^{\infty}_{j=0}$ denote an enumeration of the elements  of $I,$ listed in strictly increasing order. 
Note that $k_{j}\geq j$ for all $j.$ Further, for $j \geq 0$, set $E_{j}=A_{k_{j}, L}.$ Then
the family $\{ E_{j}\}_{ j \geq 0}$ is a partition of $G$ in finite nonempty subsets.

\smallskip
For $j \geq 0,$ set $c_{j}= c(E_{j})$, i.e.\ $ c_{j} = C_{L}(k_{j}).$ 

\medskip
We will now prove the first assertion. Using Lemma \ref{PSE}, we assume therefore that there exist $K, p > 0$ such that 
$C_{L}(k) \leq K (1+ k)^p$ for all $k \geq 0.$  Choose $s_{0} > 0$ such that  $s_{0} > p + \frac{1}{2}.$

Then we have 

\begin{align*}
\sum_{j=0}^\infty \Big(\frac{c_{j}}{(1+ k_{j})^{s_{0}}}\Big)^2 
& \leq  \sum_{j=0}^\infty K^2 \Big(\frac{(1+k_{j})^{p}}{(1+ k_{j})^{s_{0}}}\Big)^2 
= K^2  \sum_{j=0}^\infty  \frac{1}{(1+ k_{j})^{2(s_{0}-p)}} \\
& \leq K^2 \sum_{j=0}^\infty  \frac{1}{(1+ j)^{2(s_{0}-p)}} < \infty 
\end{align*}
as $2(s_{0}-p) > 1.$

Hence, defining $\kappa : G \to [1, \infty)$ by $\kappa = \sum_{j=0}^\infty (1 + k_{j})^{s_{0}} \chi_{E_{j}},$ we get from
Lemma \ref{partition} that $G$ is $\kappa$-decaying. Now, as $\kappa \leq (1+L)^{s_{0}},$
 this implies that $G$ is $(1+L)^{s_{0}}$-decaying. Assertion 1) then follows from Proposition \ref{L2D}.

\medskip
Next, assume  that $G$ has subexponential H-growth (w.r.t.\ $L$) and let $a > 1.$ 
Using Lemma \ref{PSE} we choose  $b>1$ such that $b < a$, and  $j_{0} \in \Naturali$ such that 
$C_{L}(j) < b^{\, j}$ whenever $j \geq j_{0}.$ 

\smallskip
Then we have
$$ \sum_{j=j_{0}}^\infty \Big(\frac{c_{j}}{a^{ k_{j}}}\Big)^2 
 \leq  \sum_{j=j_{0}}^\infty \Big(\frac{b^{k_{j}}}{a^{ k_{j}}}\Big)^2 
=  \sum_{j=j_{0}}^\infty  \Big(\frac{b^2}{a^2}\Big)^{k_{j}} 
 \leq  \sum_{k=k_{j_{0}}}^\infty    \Big(\frac{b^2}{a^2}\Big)^k< \infty  $$
as $b^2 / a^2 < 1.$

\smallskip
Hence, defining $\gamma : G \to [1, \infty)$ by $\gamma = \sum_{j=0}^\infty a^{k_{j}} \chi_{E_{j}},$ we get from
Lemma \ref{partition} that $G$ is $\gamma$-decaying. Now, as $\gamma \leq a^L,$
 this implies that $G$ is $a^L$-decaying. Assertion 2) then follows from Proposition \ref{L2D}.
 \hfill $\Box$

\medskip 
We may now prove Theorem \ref{I1} stated in the Introduction.
\begin{theorem} \label{M1} Let $L: G\to [0, \infty)$ be a proper function. 
   
   \smallskip If $G$ has polynomial H-growth (w.r.t.\ $L$), then there exists some $s >0 $ such that the Fourier series of 
   $x \in C_{r}^*(G, \sigma)$ 
   converges  to $x$
   in operator norm  whenever $\sum_{g\in G} | \widehat{x}(g)|^2 \, (1 + L(g))^{s} < \infty .$
   
\smallskip If $G$ has subexponential H-growth (w.r.t.\ $L$), then the Fourier series of $x \in C_{r}^*(G, \sigma)$ converges  to $x$
   in operator norm  whenever there exists some $t > 0$ such that  $\sum_{g\in G} |\widehat{x}(g)|^2 \,\exp(t L(g)) < \infty .$
   \end{theorem} 

\begin{proof} It suffices to combine Theorem \ref{Hgrowth} with Proposition \ref{L2D}, part 2).
\end{proof}

\begin{example}
\label{Hex2}
Let $G$ be any of the groups listed in Example \ref{Hex}, equipped with the length function $L$ introduced there.
As $G$ has polynomial H-growth (w.r.t.\  $L$), it follows from Theorem \ref{Hgrowth} that 
$G$ has the $\sigma$-twisted RD-property (w.r.t.\  $L$), and also that $(G, \sigma)$ is $a^L$-decaying for all $a >1$. 
\end{example}

\medskip
We conclude this section with some remarks on the interesting class of 
weight functions $\kappa$ satisfying 
$$\ \kappa(e)=1, \quad   \kappa(g^{-1})= \kappa(g), \quad \kappa(gh) \leq \kappa(g) \kappa(h)$$
for all $g,h \in G.$ Such functions are called "absolute values" in \cite{BCR}, and just "weights"
in \cite{Sch, FGL}, so we will call them \emph{absolute} weights here.
Note that $\kappa^s, s > 0$  is then also an absolute weight.
If $L$ is a length function on $G,$ then $(1+L)^s, s >0$ and $a^L, a > 1$  are all examples of such absolute weight functions.
Conversely, if $\kappa$ is an absolute weight function, then $\log_{a}(\kappa)$ is a length function for any $a>1.$

\medskip
Absolute weights are related to certain norms on $\K(G)$. If $N$ is a norm on $\K(G)$ satisfying $N(\delta_{e})=1, 
N(\xi^*)=N(\xi),$ and $N(\xi \ast_{\sigma} \eta) \leq N(\xi) N(\eta)$ for all $\xi, \eta \in \K(G),$ that is, $N$ is a 
$*$-algebra norm on $\K(G)$ (w.r.t.\   $\sigma$-twisted convolution and involution), then $\kappa_{N}(g):=N(\delta_{g})$ 
gives an absolute weight on $G.$ Conversely, one may show (using the first inequality in the next paragraph) that if $\kappa$ 
is an absolute weight on $G,$ then $N_{\kappa}:= \| \cdot \|_{1, \kappa}$ gives a norm on $\K(G)$ satisfying the above 
properties (for any $\sigma$). 

 \bigskip

Now, fix an absolute weight $\kappa$ on $G.$ For $\xi, \eta \in \L^2_{\kappa},$ it is an easy exercise to verify that
$$ |(\xi \ast_{\sigma} \eta) \kappa | \leq |\xi \kappa| \, \ast \, |\eta \kappa| .$$
This implies that 
$$ \| \xi \ast_{\sigma} \eta \|_{1, \kappa} \leq  
 \| \xi \|_{1, \kappa} \|  \eta \|_{1, \kappa}$$
whenever $\xi, \eta \in \L^1_{\kappa}.$ It follows that $\L^1_{\kappa}$ becomes a Banach $*$-algebra w.r.t.\  
$\sigma$-twisted convolution  and  involution. The problem of determining under which conditions it becomes 
symmetric as a Banach $*$-algebra has  recently been studied $G$ is of polynomial growth 
(see e.g.\ \cite{GL, FGL}). 

One may also consider $H^\infty_{\kappa}(G):= \cap_{s>0}\,\L^2_{\kappa^s},$ which becomes a 
Fr\'echet space (w.r.t.\  the obvious family of seminorms) and contains $\K(G).$
If $G$ is $\kappa$-decaying with decay constant $C,$ then we have
$$ \| \xi \ast_{\sigma} \eta \|_{2, \kappa} \leq C \| \xi \|_{2, \kappa^2} \|\eta \|_{2, \kappa}$$
whenever  $\xi \in \L^2_{\kappa^2} , \eta \in \L^2_{\kappa}$. Indeed, when
$\xi \kappa \in \L^2_{\kappa},$ we have
$$ \| \xi \ast_{\sigma} \eta \|_{2, \kappa} \leq \| |\xi \kappa| \, \ast \, |\eta \kappa| \|_{2}
\leq \| \tilde{\pi}_{\lambda}(\xi \kappa)\| \| \eta \kappa\|_{2} 
\leq C \| \xi \kappa \|_{2, \kappa} \|\eta \|_{2, \kappa}
= C \| \xi \|_{2, \kappa^2} \|\eta \|_{2, \kappa}.$$

Assume now that $G$ is $\kappa^{s_{0}}$-decaying for some $s_{0}>0.$ 
Then  one deduces from the above inequality (by considering $\xi, \eta \in H^\infty_{\kappa}(G)$
and replacing $\kappa$ with $\kappa^s$ for  $s \geq s_{0}$) 
that $H^\infty_{\kappa}(G)$ becomes a $*$-algebra under twisted convolution and involution, hence that 
 $\tilde{\pi}_{\sigma}(H^{\infty}_{\kappa}(G))$ is a (dense) $*$-subalgebra of $C^*_{r}(G,\sigma).$ 
If $\kappa = 1+L$ for some length fuction $L$ on $G$, then our assumption just says that $G$ has property RD w.r.t.\ $L,$ and 
 $\tilde{\pi}_{\sigma}(H^{\infty}_{\kappa}(G))$ is then a spectral subalgebra 
of  $C^*_{r}(G,\sigma)$ (see \cite{Cha}. and also \cite{Lu2}), as mentioned earlier in this section in the untwisted case. 
It is not unlikely that this might be generalized to more general weights.

\section{Twisted multipliers}

In \cite[Definition 1.6]{Haa1} Haagerup introduces the concept of a function which \emph{multiplies}
$C^*_{r}(G)$ into itself. The twisted analogue, which we will need in our discussion of summation processes in the next section, is as follows.

\begin{definition} 
Let $\varphi$ be a complex function on $G$.
Consider the linear map $M_\varphi: \Complessi(G, \sigma) \to \Complessi(G,\sigma)$ given by 
$$M_\varphi(\pi_{\sigma}(f)) = \pi_{\sigma}(\varphi f), \ f \in \K(G).$$
We say that $\varphi$ is a $\sigma$-\emph{multiplier}
if $M_\varphi$ is bounded w.r.t.\  the operator norm on  $\Complessi(G,\sigma)$,
in which case we also denote by $M_\varphi$ the (unique) extension of 
$M_\varphi$ to a bounded linear map from $C^*_r(G,\sigma)$ into itself.
Note that $M_\varphi$ is then uniquely determined by
$$M_\varphi(\Lambda_\sigma(g)) = \varphi(g)\Lambda_\sigma(g), \quad g \in G. $$

\noindent We denote by  $MA(G,\sigma)$ the set of all 
$\sigma$-multipliers on $G.$
Clearly $MA(G,\sigma)$ is a subspace of $\ell^{\infty}(G)$ containing
$\K(G).$ 
We set $MA(G) = MA(G,1)$, in accordance with the existing literature.
  \vspace{1ex}

\end{definition}

Adapting the arguments of Haagerup-de Canni\`ere given in the proof of \cite[Proposition 1.2]{DCHa}, 
one can show the following result.

\begin{proposition}\label{DCH}
Let $\varphi$ be a complex function on $G$. Then $\varphi  \in MA(G,\sigma)$ if and only if 
there exists a (unique) normal operator $\tilde{M}_{\varphi}$
from $vN(G,\sigma)$ to $vN(G,\sigma)$ such that
$$\tilde{M}_\varphi(\Lambda_\sigma (g)) = \varphi(g)\Lambda_\sigma(g), \quad g 
\in G, $$
in which case we have $\| M_{\varphi}\|= \| \tilde{M}_{\varphi} \|$.

Further, $MA(G,\sigma)$ is a Banach space w.r.t.\ the norm
$\tn \varphi \tn  := \| M_{\varphi} \|. $
\end{proposition}

While one implication in the first statement above is trivially true, the converse requires some work. As we won't need this result
in the sequel, we skip the details. Note however that in the course of the proof, one identifies the predual of $vN(G, \sigma)$ with a 
certain space $A(G, \sigma)$ of functions on $G$, corresponding to the Fourier algebra in the untwisted case (cf.\ \cite{Eym}), 
and establishes that $MA(G, \sigma)$ multiplies $A(G, \sigma)$ into itself. This explains the terminology and the notation.
  
\medskip Still following Haagerup-de Canni\`ere \cite{DCHa}, one
may also introduce the twisted analogue of their concept of \emph{completely bounded} multipliers :
$$M_{0}A(G,\sigma) := \{ \varphi \in MA(G, \sigma) \ | \ M_\varphi \ \text{is a completely bounded map}  \} $$
 and equip this space with the norm $\| \varphi \|_{cb}= \|M_{\varphi}\|_{cb}.$ Concerning 
completely bounded maps between C$^*$-algebras, we refer to \cite{Pau, Pis}. We set $M_{0}A(G)= M_{0}A(G, 1)$.
 
\medskip The existence of completely bounded  multipliers  
is well known in the untwisted case. Letting $P(G)$ denote the cone
of all p.d.\ functions on $G$ and $B(G) = \text{Span}(P(G))$ be the Fourier-Stieltjes algebra of $G$, then 
we have  for example $B(G) \subseteq M_{0}A(G)$ (see \cite{Haa1, DCHa, Cow, Pis}).
We recall that $B(G)$ consists of all the matrix coefficients of the unitary representations of $G$ and that 
 it may be identified with the dual space of $C^*(G)$. The norm of $\varphi \in B(G)$ as an element of 
the dual of $C^*(G)$ being denoted by $\| \varphi \|$, one has 
$\tn \varphi \tn  \leq  \| \varphi \|_{cb} \leq \| \varphi \|$. If $\varphi \in P(G),$ then
$\tn \varphi \tn  =  \| \varphi \| =\varphi(e).$ Note also that $G$ is amenable if and only if 
 $B(G)= MA(G),$ if and only if $B(G)=M_{0}A(G)$ (see \cite{N, Boz}).

 Completely bounded multipliers are closely related to (Herz-)Schur multipliers (see \cite{BF, Pis}). We recall that
a kernel $K:G \times G \to \Complessi$ is called a \emph{Schur multiplier} on $B(\ell^2(G))$ if for every
$A \in B(\ell^2(G))$ with associated matrix $[A(s,t)]$ w.r.t.\ to the canonical basis of $\ell^2(G)$, the
matrix $[K(s,t) A(s,t)]$ also represents a bounded operator on $\ell^2(G)$, denoted by $S_{K}(A)$.
When $K$ is a Schur multiplier, then the associated linear operator $S_{K}$ from $B(\ell^2(G))$ into itself
is necessarily bounded. Moreover, $S_{K}$ is then completely bounded, with $\|S_{K}\|_{cb}=\|S_{K}\|$
(see  \cite[Theorem 5.1]{Pis}).

Let now $\varphi : G \to \Complessi$ and $K_{\varphi}$ be the kernel on $G \times G$ given by 
$K_{\varphi}(s,t)= \varphi(s t^{-1})$. Then $\varphi \in M_{0}A(G)$ if and only if $K_{\varphi}$ is
a Schur multiplier, in which case we have $\|\varphi\|_{cb}=\|S_{K_{\varphi}}\|$ (see \cite{BF} and \cite[Theorem 6.4]{Pis})   In fact, 
we will show below that $\varphi \in M_{0}A(G,\sigma)$ may be characterized in the same way. Especially, we have:

\begin{proposition}\label{M0} $M_{0}A(G,\sigma)=M_{0}A(G)$ (and the cb-norm of 
$\varphi \in M_{0}A(G,\sigma)$ is independent of $\sigma$).
\end{proposition}
\begin{proof} Let $\varphi : G \to \Complessi$. As explained above, it is enough to show that $\varphi \in M_{0}A(G, \sigma)$ 
if and only if $K_{\varphi}$ is a Schur multiplier, and that in this case we have $\|\varphi\|_{cb}=\|S_{K_{\varphi}}\|$. 

Let $K_{\varphi}$ be a Schur multiplier. Then one computes that
$$[S_{K_{\varphi}}(\Lambda_{\sigma}(g))(s,t)]= [\varphi(st^{-1}) \sigma(st^{-1}, t) \delta_{g}(st^{-1})]
= [M_{\varphi}(\Lambda_{\sigma}(g))(s,t)]$$
for all $g, s, t \in G$. It follows that the restriction of $S_{K_{\varphi}}$ to $\Complessi(G, \sigma)$ is equal to
 $M_{\varphi}$. Especially, this means that $M_{\varphi}$ has a bounded extension to $C_{r}^*(G, \sigma)$, hence that
  $\varphi \in MA(G, \sigma)$. Moreover, as $S_{K_{\varphi}}$ is completely bounded,
 $M_{\varphi}$ is then also completely bounded, and 
 $\|\varphi\|_{cb}= \| M_{\varphi}\|_{cb} \leq \|S_{K_{\varphi}}\|_{cb}= \| S_{K_{\varphi}}\|.$
 
 Conversely, assume that $\varphi \in M_{0}A(G, \sigma).$ From the fundamental factorization theorem for c.b.\ maps
 (see \cite{Pau, Pis}), there exist a Hilbert space $\H$, a unital $*$-homomorphism $\pi$ from $B(\H)$ into itself, and 
 operators $T_{1}$ and $T_{2}$ from $\ell^2(G)$ into $\H$ with 
 $\|T_{1}\| \|T_{2}\| \leq \|\varphi\|_{cb}$, such that $M_{\varphi}(x)= T_{2}^* \pi(x)T_{1}$ for all
 $x \in C_{r}^*(G,\sigma)$. 
 Then 
 \begin{align*}\varphi(st^{-1})&= \varphi(st^{-1}) \, (\delta_{s}, \delta_{s})= 
 \varphi(st^{-1}) \, ( \overline{\sigma(st^{-1}, t)} \Lambda_{\sigma}(st^{-1}) \delta_{t}, \delta_{s})\\
 &=\overline{\sigma(st^{-1}, t)} \, (  M_{\varphi}(\Lambda_{\sigma}(st^{-1})) \delta_{t}, \delta_{s})
 = \overline{\sigma(st^{-1}, t)} \, (  T_{2}^*\pi(\Lambda_{\sigma}(st^{-1})) T_{1} \delta_{t}, \delta_{s})\\
 &=\overline{\sigma(st^{-1}, t)} \,  
 (  \overline{\sigma(s, t^{-1})} \pi (\Lambda_{\sigma}(s))\pi(\Lambda_{\sigma}(t^{-1})) T_{1} \delta_{t}, T_{2}\delta_{s})\\
 &= \overline{\sigma(s, t^{-1}) \sigma(st^{-1}, t) } \, (\pi(\Lambda_{\sigma}(t^{-1})) T_{1} \delta_{t}, \pi (\Lambda_{\sigma}(s))^*T_{2}\delta_{s})\\
 &= \overline{\sigma(t^{-1}, t)\sigma(s, e)} \, (\pi(\Lambda_{\sigma}(t^{-1})) T_{1} \delta_{t}, \pi (\Lambda_{\sigma}(s))^*T_{2}\delta_{s})\\
 &= (\pi(\Lambda_{\sigma}(t))^* T_{1} \delta_{t}, \pi (\Lambda_{\sigma}(s))^*T_{2}\delta_{s})
 \end{align*} 
 for all $s, t \in G$. Hence, setting $\eta_{j}(s)= \pi(\Lambda_{\sigma}(s))^* T_{j} \delta_{s} \in \H$ for $j=1,2$, we get
 $\varphi (st^{-1})= ( \eta_{1}(t), \eta_{2}(s)) $ for all $s, t \in G$, and 
 $$\sup_{s \in G}\|\eta_{1}(s)\| \, \sup_{s \in G}\|\eta_{2}(s)\|  \leq \|T_{1}\| \|T_{2}\| \leq \|\varphi\|_{cb}.$$
 From Grothendieck's theorem \cite[Theorem 5.1]{Pis}, we then deduce that 
 $K_{\varphi}$ is a Schur multiplier satisfying $\| \S_{K_{\varphi}}\| \leq \|\varphi\|_{cb}$. 
 Altogether, the assertion clearly follows. 
 \end{proof}

We don't know whether the equality  $MA(G,\sigma)= MA(G)$ always holds.

\begin{corollary}
\label{BinMA1}
We have $P(G) \subseteq B(G) \subseteq M_{0}A(G,\sigma)$. 
 If $\varphi \in P(G)$, then we have  $\tn \varphi \tn  = \| \varphi \|_{cb} = \varphi(e).$ 
\end{corollary}
\begin{proof} This follows from Proposition \ref{M0} and the facts recalled before its statement. 
As we will use the last statement several times in the sequel, we also give a direct proof.  
Let $\varphi \in P(G)$. Then write  $\varphi(\cdot) = (V(\cdot)\eta,\eta)$ for some 
unitary representation  $V$ of $G$ on a Hilbert space $\H$ and  some $\eta \in \H$.
Let $W$ be the unitary operator on $\ell^2(G) \otimes \H \cong \ell^2(G,\H)$ given by  
$$( W \psi)(g) = V(g) \psi(g), 
\ g \in G, \psi \in \ell^2(G, \H).$$ 
Then one computes that 
$W^*(\Lambda_\sigma (g) \otimes V(g) )W =  \Lambda_{\sigma}(g) \otimes I_\H $
for all $g \in G$. 
This is the twisted version (cf.\ \cite[Prop. 2.2]{BeCo}) of usual Fell's absorbing property (\cite[13.1.3]{Di}).

Next, let $T: \ell^2(G) \to \ell^2(G) \otimes \H$ be  given by 
$T(\xi) = \xi \otimes \eta$. 
Then  $T$ is linear, bounded  and 
$T^*(\xi' \otimes \eta') = (\eta',\eta) \xi'.$ 

Now, define $M: C_{r}^*(G,\sigma) \to B(\ell^2(G))$ by 
$$M(x) = T^* W (x \otimes I_{\H}) W^* T.$$
Then $M$ is  a completely positive map (see \cite{Pau}) and 
$$(*) \quad\|M\| = \|M \|_{cb}= \|M(I) \| = \|T^*T\| = \|\eta\|^2   =  \varphi (e).$$
Furthermore, for $g \in G, \xi \in \ell^2(G)$, we have
\begin{align*}
M(\Lambda_{\sigma}(g))\xi & = T^*W(\Lambda_{\sigma}(g) \otimes I_{\H})W^*T \xi 
 = T^* (\Lambda_\sigma(g)\xi \otimes V(g)\eta) \\
& = (V(g)\eta,\eta) \Lambda_\sigma(g)\xi  = \varphi(g) \Lambda_\sigma(g) \xi 
 \end{align*}
Hence, it follows that $M$ is a c.b.\ extension of $M_{\varphi}$ and the last statement follows from $(*)$. 
\end{proof}

\begin{remark} $\ell^2(G) \subseteq M_{0}A(G,\sigma)$ (with $\| \varphi \|_{cb} \leq \| \varphi \|_{2}$).
This is easy to see directly, but also follows from Proposition \ref{M0} (as $\ell^2(G) \subseteq B(G)$).
\end{remark}

Now, to prepare for our study of summation processes in the next section, 
consider   $\varphi \in MA(G,\sigma)$  and 
$x \in C^*_r(G,\sigma)$. Then $\widehat{M_\varphi(x)} = \varphi \widehat{x}$.  
 
\smallskip \noindent Indeed, if $x \in \Complessi(G,\sigma),$ this is trivial; 
otherwise the statement follows immediately from a density argument. Hence, the Fourier series of $M_\varphi(x)$ is 
$\sum_{g\in G} \varphi(g) \widehat{x}(g)\Lambda_\sigma(g)$.
This series does not necessarily converge  in operator norm, but if for example $\varphi \in \ell^2(G)$,  
then it does, since $\varphi \widehat{x} \in \ell^1(G)$. This motivates the following definition.

\begin{definition} We let 
$MCF(G,Ê\sigma)$ denote the set of all  complex functions $\varphi : G \to \Complessi $
such that the series $\sum_{g\in G} \varphi(g)\widehat{x}(g)\Lambda_{\sigma}(g)$
converges in operator norm for all $x \in C^*_{r}(G, \sigma).$
\end{definition}
At least, we know that $\ell^2(G) \subseteq MCF(G,Ê\sigma)$. Further, we have the following.

\begin{proposition} \label{MCF} 
 $MCF(G,Ê\sigma) \subseteq MA(G, \sigma).$ Moreover, 
$$MCF(G,Ê\sigma) = \{ \varphi \in MA(G, \sigma)  \, |  \,  M_{\varphi}  \
 \text{maps}\ C^*_{r}(G,\sigma) \ \text{into} \ CF(G,\sigma) \}$$
 and if $\varphi \in MCF(G,Ê\sigma)$, then $ \sum_{g\in G} \varphi(g)\widehat{x}(g)\Lambda_{\sigma}(g) $ 
 converges to $M_{\varphi}(x)$ in operator norm for  all $ x \in C^*_{r}(G, \sigma).$
\end{proposition}
\begin{proof} Let $\varphi \in MCF(G,Ê\sigma).$ Define a linear map $M'_{\varphi}: C^*_{r}(G, \sigma)
 \to C^*_{r}(G, \sigma)$ by 
 $$M'_{\varphi}(x) = \sum_{g\in G} \varphi(g)\widehat{x}(g)\Lambda_{\sigma}(g).$$
Using the closed graph theorem, one gets that  $M'_{\varphi}$ is bounded.  
Indeed, assume $x_{n} \to x$  and $M'_{\varphi}(x_{n}) \to y$ in $C^*_{r}(G,\sigma).$ Then 
$\widehat{M'_{\varphi}(x_{n})}=\varphi \widehat{x_{n}} \to  \varphi \widehat{x}= \widehat{M'_{\varphi}(x)}$
 pointwise on $G$ and also $\widehat{M'_{\varphi}(x_{n})} \to \widehat{y}$ pointwise on $G.$ 
 Hence, $\widehat{M'_{\varphi}(x)}= \widehat{y},$ so $M'_{\varphi}(x)= y$, as desired.

As $M'_\varphi(\Lambda_\sigma(g)) = \varphi(g) \Lambda_\sigma(g)$ for all $ g \in G,$ this implies that $M'_{\varphi}$ 
is a bounded extension of  $M_{\varphi}$ from $\Complessi(G, \sigma)$ to  $C^*_{r}(G, \sigma).$ 
 Hence $\varphi \in MA(G,Ê\sigma)$, and the first statement is proven.  
 As $\widehat{M_\varphi(x)} = \varphi \widehat{x}$ for all $ x \in C^*_{r}(G, \sigma)$, the last assertion follows.
 \end{proof}
 
Inspired by \cite[Lemma 1.7]{Haa1}, we can produce other examples of multipliers 
in $MCF(G, \sigma)$.

\begin{proposition}
\label{v1.7}
Let $G$ be $\kappa$-decaying with decay constant $C.$ \\
\smallskip Let $\psi \in \L^{\infty}_{\kappa}$ and set $K = \|\psi\|_{\infty, \kappa}$. 
Then $\psi \in MCF(G, \sigma)$ with $\tn  \psi \tn  \leq C K$.
\end{proposition}

\begin{proof} From Proposition \ref{L2D} we know that $(G, \sigma)$ is $\kappa$-decaying. Moreover, 
from the proof of Proposition \ref{norm}, we see that $\|\pi_{\sigma}( f)\| \leq C \| f\|_{2,\kappa}$
for all $f \in \K(G)$,  where $C$ is given as above.

\smallskip Now, let $f \in \K(G).$ Then
$$ \|\psi f\|_{2,\kappa}
=  \| \psi f \kappa \|_{2} \leq  \|\psi \kappa\|_\infty \|f\|_2  = K \|f\|_2  \leq K \|\pi_{\sigma}(f)\|.$$
 Hence we get
$$\|M_\psi(\pi_{\sigma}(f))\| = \|\pi_{\sigma}(\psi f)\| \leq C \|\psi f\|_{2,\kappa}
\leq CK \|\pi_{\sigma}(f)\|.$$
Thus $M_{\psi}$ is bounded with $\|M_\psi\| \leq C K.$ Especially, $\psi \in MA(G, \sigma)$ and 
it remains only to show that $\psi \in MCF(G, \sigma).$

\smallskip Let $x \in C^*_{r}(G, \sigma).$ As
$\| \psi \widehat{x}\|_{2,\kappa} \leq  K \|\widehat{x}\|_{2} < \infty$, we have
$\widehat{M_{\psi}(x)} = \psi \widehat{x} \in \L^2_{\kappa}$. 
From the last statement in Proposition \ref{L2D}, we get that 
$\sum_{g\in G} \psi(g)\widehat{x}(g)\Lambda_{\sigma}(g) $ converges in operator norm, as desired. 
\end{proof}

\medskip
\begin{remark} \label{Feller}
Let $G$ be $\kappa$-decaying and let $\psi \in \L^{\infty}_{\kappa}$. The proof of
Proposition \ref{v1.7}  shows in fact that $\sum_{g \in G}\psi(g)\widehat{x}(g)\Lambda_{\sigma}(g)$
is operator norm convergent for all $x \in vN(G, \sigma).$  

If, in addition, $\psi$ is  p.d., 
then it is natural to wonder whether it has the strong Feller property 
introduced by J.L.\ Sauvageot \cite{Sau, Sau2}, that is, whether  
$M_{\psi}^{**} (C_{r}^*(G, \sigma)^{**}) \subseteq C_{r}^*(G, \sigma).$ 
Now, one readily sees from the proof of Proposition \ref{v1.7} that there exists a constant $C' >0$ such that 
 $\|\pi_{\sigma}(\psi f)\| \leq C' \| fÊ\|_{2} $ for all $f \in \K(G),$ and it does indeed follow
  that $\psi$ has the strong Feller property (cf.\ \cite[Lemma 3.3 and Proposition 5.2]{Sau2}).
\end{remark}

\section{Summation processes}
We begin with some definitions. 

\smallskip \begin{definition}
 A net $\{\varphi_\alpha\}$  in $MA(G, \sigma)$
is called an \emph{approximate multiplier unit} 
whenever $M_{\varphi_\alpha}(x) \to x$ in operator norm for every $x \in C^*_r(G,\sigma)$.

\noindent Such a net is called \emph{bounded}  when $\sup_{\alpha} \tn \varphi_{\alpha} \tn < \infty $.
\end{definition}

\begin{remark}
\label{remamu}
We record the following simple but useful facts :
\begin{itemize}
\item[1)] 
 Assume that $\{\varphi_\alpha\}$  is an approximate multiplier unit in $MA(G, \sigma)$.
Then $\varphi_\alpha \to 1$ pointwise on $G$ and we have 
$1 \leq \sup_\alpha \tn \varphi_\alpha \tn \leq \infty.$ 
If $\{\varphi_\alpha\}$ is a sequence, then $\{\varphi_\alpha \}$ is bounded  
(as follows from the uniform boundedness principle).
\item[2)] 
Let $\{\varphi_\alpha \}$ be a net in $MA(G,\sigma).$ 
Using a straightforward $\varepsilon /3$-argument, one deduces that 
$\{\varphi_\alpha\}$  is a bounded approximate multiplier unit if and only if
 $\varphi_\alpha \to 1$ pointwise on $G$ and $\{\varphi_\alpha\}$ is bounded.
\item[3)]   
  If $G$ is countable (so $C_{r}^*(G,\sigma)$ is separable), then (mimicking the trick used to produce
  a countable approximate unit in a separable $C^*$-algebra) one can always extract a sequence from a given bounded
  approximate multiplier unit to produce a (bounded) countable approximate multiplier unit if necessary.

\end{itemize}
\end{remark}

\begin{example} \label{examu}
Assume that 
$\{ \varphi_{\alpha}\}$ is a net of normalized p.d.\ functions on $G$ converging pointwise to $1.$
Then $\tn \varphi_{\alpha}\tn =1$ for all $\alpha$ (cf.\ Corollary \ref{BinMA1})
and  assertion 2) above gives that $\{ \varphi_{\alpha}\}$ is a bounded approximate multiplier unit for $C_{r}(G,\sigma).$
\end{example}

\begin{definition}
Let  $\{\varphi_\alpha\}$ be a net of complex functions on $G.$
We say that $\{\varphi_\alpha\}$ is a \emph{Fourier summing net for} $(G, \sigma)$
 if $\{\varphi_\alpha\}$ is an approximate multiplier unit for $C_{r}(G,\sigma)$ satisfying
 $\varphi_{\alpha} \in MCF(G,\sigma)$ for  all $\alpha.$

\smallskip Such a net gives a summation process for Fourier series of elements in $C^*_{r}(G,\sigma)$ :  
the series $\sum_{g \in G} \varphi_\alpha(g) \widehat{x}(g) \Lambda_\sigma(g)$ is then
convergent in operator norm for all $\alpha,$ and 
 $$\sum_{g \in G} \varphi_\alpha(g) \widehat{x}(g) \Lambda_\sigma(g) \underset\alpha\to x$$ 
for all $x \in C^*_r(G,\sigma)$ (w.r.t.\ operator norm). 
\end{definition}

It is an open question whether one can always find a Fourier summing net for a general pair $(G,\sigma).$ 
When $G$ is amenable, the answer is well-known. Indeed, the following theorem was proven  by Zeller-Meier in \cite{ZM} (see also \cite{Exe}) 
in the case of a net of finitely supported functions.

\begin{theorem}\label{ZM}
Let $G$ be amenable and $\{\varphi_\alpha\}$ be  any net of  normalized p.d.\ 
functions in $\ell^2(G)$  converging pointwise to $1$. 
Then  $\{\varphi_\alpha\}$ is a  (bounded) Fourier summing net for $(G,\sigma)$ 
(satisfying  $\tn \varphi_\alpha \tn = 1$ for all $\alpha$).
\end{theorem}

\begin{proof}  As $\ell^2(G) \subseteq MCF(G,\sigma)$ (cf.\ Section 4), this follows from Example \ref{examu}. 
\end{proof}

\medskip We turn now to the proof of Theorem \ref{I2} on Fej\'er summation, which may be restated as follows : 

\begin{theorem}\label{Fej} Let $G$ be amenable and pick a F\o lner net  $\{F_\alpha\}$ for $G.$  Set 
$$\varphi_\alpha(g) = \frac{|g F_\alpha  \cap F_\alpha|}{|F_\alpha|}, 
\quad  g \in G.$$ (Note that each $\varphi_\alpha$ has finite support given by 
${\rm supp}(\varphi_\alpha) = F_\alpha \cdot F_\alpha^{-1}$). Then 
$\{\varphi_\alpha\}$ is a  (bounded) Fourier summing net for $(G,\sigma)$ (satisfying  $\tn \varphi_\alpha \tn = 1$ for all $\alpha$).
\end{theorem}
 \begin{proof} Each $\varphi_\alpha$ is normalized, and
the F\o lner condition gives that $\varphi_\alpha$ converges pointwise to $1.$ 
As  $\varphi_{\alpha}(g)= (\lambda(g)\xi_{\alpha},\xi_{\alpha} ),$ where
 $\xi_{\alpha}:= |F_{\alpha}|^{-1/2} \chi_{F_{\alpha}},$ 
each $\varphi_{\alpha}$ is positive definite. This means that
 $\{\varphi_\alpha\}$ satisfies the assumptions of Theorem \ref{ZM} and the result follows.
\end{proof} 
We remark that N.\ Weaver \cite{Wea} has proved this result for twisted group $C^*$-algebras of $\Relativi^2$, 
using a different approach.   
 
  \medskip Next, we turn our attention to Abel-Poisson summation and prove first Theorem \ref{I3}. We restate it in  a 
  slightly more general form, which also incorporates Gauss summation. 
  
   \begin{theorem}  Let $ G=\Relativi^N$ for some $N \in \Naturali$. For $p \in \{1, 2\}$,
    let $| \cdot |_{p}$ denote the usual $p$-norm on $G$. Let $L(\cdot)$ denote either $| \cdot |_{1}\, , | \cdot |_{2}$ or 
    $| \cdot |_{2}^2$. For each $r \in (0,1)$, set $\varphi_{r} = r^{L}$. 
   
   Then $\{\varphi_{r}\}_{r \to 1^{-}}$ is a (bounded) Fourier summing net for $(G,\sigma)$. 
\end{theorem} 

 \begin{proof} It is well known and elementary that $| \cdot |_{2}^2$ is n.d.\ . Hence,  $| \cdot |_{2}$, being the square root
 of $| \cdot |_{2}^2$, is also n.d.\ (see \cite{BCR}). Especially, $| \cdot |_{1}$ is n.d.\ when $N=1$, and it follows from
 a simple inductive argument that $| \cdot |_{1}$ is n.d.\ for all  $N\geq 1$. This means that  $L$ is n.d.\  Hence, according to Theorem \ref{S},
 all $\varphi_{r}$ are
 p.d.\  Morover, one checks easily that they are square-summable and normalized. As 
$\varphi_{r}$  converges pointwise to 1 when $r \to 1^{-}$, Theorem \ref{ZM} applies and gives the result.
   \end{proof}
   
   The Gaussian case above (choosing $L(\cdot)=| \cdot |^2_{2}$) illustrates that one should not only consider length functions. 
To show the existence of summation processes for many other (nonamenable) groups, we will use the following result. 

\begin{proposition}
\label{remsp}
 Let $\{\varphi_\alpha\}$ be a net in $MA(G, \sigma)$. Assume that 
  \begin{itemize}
\item[i)] $\{\varphi_\alpha\}$ converges pointwise to $1$,  
 \item[ii)] $\{\varphi_\alpha\}$ is bounded,
\item[iii)]  for each $\alpha$ there exists some $\kappa_{\alpha} : G \to [1, \infty)$ such that 
 $G$ is $\kappa_{\alpha}$-decaying
 and $\{\varphi_\alpha\} \in \L^{\infty}_{\kappa_{\alpha}}.$
\end{itemize}
\noindent 
\smallskip Then $\{\varphi_\alpha\}$ is a (bounded) Fourier summing net for $(G, \sigma)$.
\end{proposition}

\begin{proof} Conditions i) and ii) ensure that $\{\varphi_\alpha\}$ is a bounded approximate multiplier unit (cf.\ Remark \ref{remamu}, part 2)).
Further, Proposition \ref{v1.7} ensures that $\{\varphi_\alpha\} \subseteq MCF(G,\sigma).$ 
\end{proof}

\begin{theorem} \label{APgrowth} Let $G$ be a countable group with the Haagerup property and 
$L$ be a Haagerup function on $G$. 

 \smallskip \noindent
 Assume that $G$ has polynomial H-growth (w.r.t.\  $L$). Then there exists some $q \in \Naturali$ 
 such that $\{ (1 + t L )^{-q} \}_{t \to 0^{+}}$ is a (bounded) Fourier summing net for  $(G, \sigma)$.   

\smallskip \noindent
More generally, assume that $G$ has subexponential H-growth (w.r.t.\  $L$). 
Then  $\{r^{L}\}_{r \to 1^{-}}$ is  a (bounded) Fourier summing net for  $(G, \sigma)$.
\end{theorem}

\begin{proof} For $ p \in \Naturali, t >0, $ set   $\psi_{p,t}= (1 + t L)^{-p}$. 
  For $r \in (0,1)$, set $\varphi_{r}= r^L$. Then $\psi_{p,t}$   and $\varphi_{r}$  are normalized positive
definite functions on $G$, as follows respectively from \cite[p.\ 75]{BCR} and from Theorem \ref{S}. Hence,  both 
$\{ \psi_{p,t} \}_{t \to 0^{+}}$ and $\{\varphi \}_{r \to 1^{-}}$ are bounded  (cf.\ Example \ref{examu}) and  
converge pointwise to 1.

 \smallskip
Assume first that $G$ has polynomial H-growth w.r.t.\  $L$. Due to   Theorem \ref{Hgrowth}, part 1), we may 
pick $s_{0}>0$  such that $G$ is $\kappa$-decaying, where $\kappa = (1+L)^{s_{0}}.$ 
Choose $q \in \Naturali$ such that $q \geq s_{0}.$ Clearly, $\psi_{q,t} \in \L^\infty_{\kappa} $  for all $t >0$.
This means that $\{\psi_{q,t}\}_{t\to 0^{+}}$ satisfies all conditions in Proposition  \ref{remsp} 
(with $\kappa_{t}= \kappa$ for all $t>0$), and the first assertion follows.

\smallskip
Next, assume that $G$ has subexponential H-growth (w.r.t.\   $L$). 
Let $0<r<1$, set $\kappa_{r}=r^{-L}$.  Then, according to Theorem \ref{Hgrowth}, part 2),
$(G, \sigma)$ is $\kappa_{r}$-decaying. 
Moreover, we obviously have $\varphi_{r} \in \L^\infty_{\kappa_{r}}.$   
This means that  $\{\varphi_{r}\}_{r \to 1^{-}}$ satisfies all conditions in Proposition  \ref{remsp}, 
and the second assertion follows.
\end{proof}

\begin{example} \label{niceg}
Let $G$ be a finitely generated free group, or a Coxeter group, with generator set $S$.
Then the word-length $L_{S}$ is a Haagerup function on $G$ (see \cite{CCJJV}).
Further,  $G$ has polynomial H-growth w.r.t.\ $L_{S}$ (see Example 
\ref{Hex}). Hence,  Theorem \ref{APgrowth} applies.
\end{example}

\begin{remark} Assume that there exists a net $\varphi_{\alpha}$ of normalized p.d.\ functions on 
 $G$ converging pointwise to $1$ and satisfying condition iii)  in Proposition \ref{remsp}.
Then $C_{r}^*(G,\sigma)$ has the strong Feller approximation property considered by Sauvageot \cite{Sau, Sau2}: indeed, each
$\varphi_{\alpha}$ has then the strong Feller property, cf.\ Remark \ref{Feller}. 
This observation applies to any countable group which has the Haagerup property
 and has  subexponential H-growth w.r.t.\ some Haagerup function 
 (cf.\ Theorem \ref{APgrowth} and its proof). 
 \end{remark}

The class of groups for which the Abel-Poisson summation holds contains indeed many other groups.
 
\begin{theorem} \label{Hyp}
Let $G$ be a Gromov hyperbolic group and let $L$ be an algebraic length function on $G$.
Then  $\{r^{L}\}_{r \to 1^{-}}$ is  a (bounded) Fourier summing net for  $(G, \sigma)$.
\end{theorem}
\begin{proof} In a recent paper \cite{Oz}, N.\ Ozawa has shown that the net  $\{r^{L}\}_{r \to 1^{-}}$ is c.b.\ 
 bounded in $M_{0}A(G)$. Using Proposition \ref{M0}, we get that this net is c.b.\ bounded in $M_{0}A(G,\sigma)$. 
 In particular, it is bounded in $MA(G, \sigma)$. Moreover, as explained in Example \ref{Hex}, $G$ has polynomial 
 H-growth, hence subexponential H-growth (w.r.t.\ $L$). We can now conclude the proof by proceeding in the same way as 
 in the proof of the second statement of Theorem \ref{APgrowth}.  
\end{proof}

We conclude this paper with some remarks on Fej\'er-like properties.

\begin{definition}
We say that $(G,\sigma)$ has the {\emph Fej\'er property} 
if there exists a  Fourier summing net $\{\varphi_\alpha\}$ for $(G, \sigma)$ in $\K(G)$.
If  $\{\varphi_\alpha\}$ converges pointwise to $1$ and is bounded in $MA(G,\sigma)$,
then we  say that $(G,\sigma)$ has the {\emph bounded Fej\'er property}.
Moreover, if this net can be chosen to satisfy $\tn \varphi_\alpha \tn = 1$ for all $\alpha$, then we  say that
$(G,\sigma)$ has the {\emph metric Fej\'er property}.

When $\sigma =1,$ we just talk about the corresponding Fej\'er property for the group $G.$
\end{definition}

 To motivate the use of the adjective ``metric'' in the metric Fej\'er property,
we recall that a Banach space $X$ is said to have the \emph{Metric
Approximation Property} (M.A.P.) if there exists a net of finite rank contractions on $X$ 
approximating the identity map in the strong operator topology (SOT) on $B(X).$ 
Hence, if  $(G,\sigma)$  has the metric Fej\'er property, 
then $C^*_{r}(G,\sigma)$ has the M.A.P.\ We don't know whether the converse is true.
In  \cite[Theorem 1.8]{Haa1}, Haagerup shows that $\mathbb{F}_2$ has the metric Fej\'er property, 
hence that $C^*_r({\mathbb F}_2)$ has  the M.A.P. (despite the fact that $C^*_r({\mathbb F}_2)$ is not nuclear).

Theorem \ref{Fej} shows that $(G,\sigma)$ has the metric Fej\'er property
whenever $G$ is amenable. It is not unlikely that this is still true whenever $G$ has the Haagerup property. In fact, 
in the untwisted case, a conjecture of M.\ Cowling (see \cite{CCJJV}) says that any (countable) group $G$ 
with the Haagerup property is \emph{weakly amenable} \cite{CH} with CH-contant equal to 1, that is, there exists a
net $\{\varphi_\alpha \}$ in $\K(G)$ converging pointwise to $1$ such that 
 $\sup_\alpha\| \varphi_\alpha \|_{cb} = 1$. Cowling's conjecture (which also may be formulated
 in the locally compact case) has been verified in a number of cases (see e.g.\ \cite{CH, Cow2, Va2, Jan, CCJJV, GuHi2}). 
 
The following result generalizes \cite[Theorem 1.8]{Haa1} (see also \cite{JoVa} and \cite{BrNi}). 
\begin{theorem}
\label{genHaagerup}
Assume that the following  conditions hold:
\begin{itemize}
\item[(i)] 
There exists a net $\{\varphi_\alpha\}$ in $MA(G,\sigma)$ converging pointwise to 1 and
satisfying $\tn \varphi_\alpha \tn = 1$ for all $\alpha$. 
\item[(ii)] 
For each $\alpha$ there exists a function $\kappa_{\alpha}: G \to [1, + \infty)$ such that  $G$ is $\kappa_{\alpha}$-decaying
 and  $\varphi_\alpha \kappa_{\alpha} \in c_0(G)$.
\end{itemize}
\noindent Then $(G,\sigma)$ has the metric Fej\'er property. 
\end{theorem}

\begin{proof} 
Clearly, $\varphi_\alpha \neq 0$ for all $\alpha$. Let $\alpha \in \Lambda$, $n \in \Naturali$.
Using (ii), we can pick a finite subset $A_{\alpha,n}$ of $G$ such that 
$|\varphi_\alpha \kappa_{\alpha}| \leq \frac{1}{n}$ outside $A_{\alpha,n}$.
If necessary, we enlarge  $A_{\alpha,n}$ to include at least one element where $\varphi_\alpha$ is nonzero.
Set $\varphi_{\alpha,n} = \varphi_\alpha \chi_{A_{\alpha,n}} $. Then 
$$\|(\varphi_\alpha - \varphi_{\alpha,n}) \kappa_{\alpha} \|_\infty
= \sup\{|(\varphi_\alpha \kappa_{\alpha})(g)| \ , \ g \notin A_{\alpha,n}\} 
\leq \frac{1}{n} \ .$$
Using  Proposition \ref{v1.7}, we get that $(\varphi_\alpha - \varphi_{\alpha,n}) \in MA(G, \sigma)$ and
$$\tn \varphi_\alpha - \varphi_{\alpha,n} \tn \leq \frac{C_{\alpha}}{n}\to 0 \ \text{as} \ n \to +\infty $$ 
where $C_{\alpha}$ denotes the decay constant of $G$ w.r.t.\  $\kappa_{\alpha}.$

\smallskip
So, setting $\psi_{\alpha,n}  = \frac{1}{\tn \varphi_{\alpha,n} \tn} \varphi_{\alpha,n},$
we have $\tn \psi_{\alpha,n} \tn = 1$ and $\tn \psi_{\alpha,n} - \varphi_\alpha \tn \to 0$ as $n \to \infty$.
 
 Now, using (i), we have $M_{\varphi_\alpha} \to {\rm Id}$ in the SOT on $B(C^*_r(G,\sigma))$. 
It follows easily that ${\rm Id} \in \{M_{\psi_{\alpha,n}} \, | \, \alpha \in \Lambda, n \in \Naturali\}^{- {\rm SOT}}.$ 
The existence of a net $\{\psi_\beta\}$ in  $\K (G)$ converging pointwise to $1$ and  
satisfying  $\tn \psi_\beta \tn = 1$ for all $\beta$ is then clear. Hence, $(G, \sigma)$ has the metric Fej\'er property
\end{proof}

\begin{corollary}\label{last}
Assume that $G$ is countable and has the Haagerup property. If  there exists  a Haagerup function $L$ on $G$ such that
$G$ has subexponential H-growth (w.r.t.\   $L$),  then $(G,\sigma)$ has the metric Fej\'er property.
  \end{corollary}

\begin{proof}
Assumptions (i) and (ii) in Theorem \ref{genHaagerup} hold with 
$\varphi_r= r^{L}$ and $\kappa_{r}= (r^{-1/2})^L$, $ 0 < r \to 1^{-}$, and the result follows.
\end{proof}

 N.\ Ozawa has recently shown \cite{Oz} that all Gromov hyperbolic groups are weakly amenable, hence especially they have the bounded 
 Fej\'er property. In a certain sense, "most" finitely presented groups are Gromov hyperbolic (see \cite{O}). However, 
 not all groups have the bounded Fej\'er property. This follows from an unpublished work of Haagerup \cite{Haa3}, where he considers 
the group $H$ obtained by forming the standard semi-direct product  of $\Relativi^2$
by  $SL(2, \Relativi)$ and shows that $H$ is not weakly amenable by actually proving that $H$ does not have the
bounded Fej\'er property (see \cite{Do} for the continuous version of this result). But note that $H$, which does not have the
Haagerup property, does have the Fej\'er property: this follows from \cite{HK}, where Haagerup and Kraus  show that $H$ 
satisfies a certain approximation property, called AP, which is stronger than the Fej\'er property. It is conceivable that  $SL(3, \Relativi)$ 
does not have the Fej\'er property. Haagerup and Kraus conjecture in \cite{HK} that $SL(3, \Relativi)$
fails to have the AP, but this is still open as far as we know. 

\smallskip
\noindent{\bf Acknowledgements.} 
This work started while the second author was visiting the University of Oslo in April 2005, in connection with the SUPREMA-program.
It was continued during a second visit, in April/May 2006. He wishes to thank all the members of the Operator Algebras group in 
Oslo for the kind invitations and hospitality. 

Both authors would also like to thank the referees for their valuable comments after reading 
 previous drafts of this paper, and for their suggestions concerning the literature on this subject.

\end{document}